\documentclass[twoside,11pt,reqno]{amsart}
\usepackage{amsmath,amssymb,amscd,mathrsfs}

\def\@underbar#1#2{\settowidth{\@tempdimb}{$#1#2$}\@tempdimb=0.8\@tempdimb
                   \ooalign{$#1#2$\crcr
                         \hfil\rule[-.5mm]{\@tempdimb}{.4pt}\hfil}}
\newcommand{\Fr}{F}
\newcommand{\wt}{\operatorname{wt}}
\newcommand{\ui}{{\text{\boldmath{$i$}}}}
\newcommand{\uj}{{\text{\boldmath{$j$}}}}
\newcommand{\uk}{{\text{\boldmath{$k$}}}}
\newcommand{\ul}{{\text{\boldmath{$l$}}}}
\newcommand{\us}{{\text{\boldmath{$s$}}}}
\newcommand{\ut}{{\text{\boldmath{$t$}}}}
\newcommand{\uh}{{\text{\boldmath{$h$}}}}
\newcommand{\p}[1]{\ensuremath{\overline{#1}}}
\newcommand{\losemi}{{\otimes \kern -.78em \ltimes}}
\newcommand{\rosemi}{{\otimes \kern -.78em \rtimes}}

\newcommand{\Hom}{\ensuremath{\operatorname{Hom}}}
\newcommand{\End}{\ensuremath{\operatorname{End}}}
\newcommand{\Coind}{\ensuremath{\operatorname{coind}}}

\newcommand{\Res}{\ensuremath{\operatorname{res}}}

\newcommand{\Ker}{\ensuremath{\operatorname{Ker} }}

\newcommand{\I}{\ensuremath{\mathcal{I}_{1}}}
\newcommand{\dist}[1]{\operatorname{Dist}_{#1}(G)}
\newcommand{\Dist}[1]{\operatorname{Dist}(#1)}

\newcommand{\distev}[1]{\operatorname{Dist}_{#1}(G_\ev)}
\newcommand{\lie}{\operatorname{Lie}(G)}

\newcommand{\sgn}{\operatorname{sgn}}
\newcommand{\0}{\bar 0}
\newcommand{\1}{\bar 1}

\newcommand{\mx}{\ensuremath{\mathtt{M}}}
\newcommand{\rx}{\ensuremath{\mathtt{R}}}

\newcommand{\jx}{\ensuremath{\mathtt{J}}}

\newcommand{\gl}{\ensuremath{\mathfrak{gl}(m|n)}}
\newcommand{\ev}{{\operatorname{ev}}}
\newcommand{\Gev}{G_{\ev}}
\newcommand{\Gevr}{\ensuremath{G_{\text{ev},r}}}
\newtheorem{Df}{Def.}
\newtheorem{theorem}[Df]{Theorem}

\newtheorem{lemma}[Df]{Lemma}

\newtheorem{remark}[Df]{Remark}

\newtheorem{corollary}[Df]{Corollary}
\numberwithin{Df}{section}
\numberwithin{equation}{section}

\begin{document}
\title{A new proof of the Mullineux conjecture}
\author{Jonathan Brundan and Jonathan Kujawa}
\address{Dept. of Mathematics \\
           University of Oregon \\
           Eugene, OR 97401}
\email{brundan@darkwing.uoregon.edu, kujawa@noether.uoregon.edu}
\thanks{Research supported in part by the NSF (grant no DMS-0139019)}

\maketitle

\section{Introduction}

Let $S_n$ be the symmetric group on $n$ letters, $k$ be a field of characteristic $p$ and $D^\lambda$
be the irreducible $kS_n$-module 
corresponding to a $p$-regular partition $\lambda$ of $n$, as in \cite{J2}.
By tensoring $D^\lambda$ with the $1$-dimensional sign representation
we obtain another irreducible $kS_n$-module.
If $p = 0$, $D^{\lambda} \otimes \sgn \cong D^{\lambda'}$,
where $\lambda'$ is the conjugate of the partition $\lambda$, and
if $p = 2$, we obviously have that
$D^{\lambda} \otimes \sgn \cong D^{\lambda}$. 
In all other cases, it is surprisingly difficult to describe
the partition labeling the irreducible module
$D^{\lambda} \otimes \sgn$ combinatorially.
In 1979, Mullineux \cite{M1} gave an algorithmic construction of a bijection
$\mx$ on $p$-regular partitions, 
and conjectured that $D^{\lambda} \otimes \sgn \cong D^{\mx(\lambda)}$.

Mullineux's conjecture was finally proved in 1996.
The key breakthrough leading to the proof was made in \cite{K2},
when Kleshchev discovered an alternative algorithm,
quite different in nature to Mullineux's,
and proved that it computes the label of $D^{\lambda} \otimes \sgn$.
Then Ford and Kleshchev \cite{FK1} proved combinatorially
that Kleshchev's algorithm was equivalent to Mullineux's, hence
proving the Mullineux conjecture.
Since then, different and easier approaches to 
the combinatorial part of the proof, i.e. that Kleshchev's
algorithm equals Mullineux's algorithm,
have been found by Bessenrodt and Olsson \cite{BO2} and by Xu \cite{Xu2}.
Also  Lascoux, Leclerc and Thibon \cite{LLT} 
have used Ariki's theorem \cite{A} to give
a different proof of the results of \cite{K2}.

The purpose of the present article is to explain a completely 
different proof of the Mullineux conjecture.
In \cite{Xu1}, Xu discovered yet another algorithm,
and gave a short combinatorial argument to show that it
was equivalent to Mullineux's original algorithm.
We will show directly from representation theory 
that Xu's algorithm computes the label of $D^{\lambda} \otimes \sgn$.
In this way, we obtain a relatively direct proof of
the Mullineux conjecture that bypasses Kleshchev's algorithm
altogether.

The idea behind our approach is a simple one. There is a superalgebra analogue
of Schur-Weyl duality relating representations of $S_n$
to representations of the supergroup $GL(n|n)$. 
Moreover, 
there is an involution on representations of $GL(n|n)$ 
induced by twisting with its 
natural outer automorphism, which corresponds under
Schur-Weyl duality to tensoring with the sign representation.
Ideas of Serganova \cite{sergthesis}
give an easy-to-prove algorithm for
computing this involution, hence by Schur-Weyl duality
we obtain an algorithm for computing $D^{\lambda} \otimes \sgn$.
Actually, we obtain a whole family of algorithms, 
one of which turns out to be the same as Xu's algorithm.

The remainder of the article is arranged as follows.
In $\S$\ref{S:basics}, we review some generalities concerning the
supergroup $G = GL(m|n)$. In $\S$\ref{S:algebraofdistributions}, 
we introduce the superalgebra $\Dist{G}$ of distributions on $G$ 
and explain how integrable representations of $\Dist{G}$ can be lifted
to $G$ itself.
Serganova's algorithm is derived in $\S$\ref{S:serg} using some
highest weight theory.
In $\S$\ref{S:schuralgebra}, we review some known results about polynomial
representations and Schur-Weyl duality,
allowing us to descend to the symmetric group. 
Finally in $\S$\ref{S:Mullineux} 
we put it all together with some combinatorics
to obtain the proof of the Mullineux conjecture.
At the end of $\S$\ref{S:Mullineux}, we also solve a related
question concerning the classification of 
the irreducible polynomial representations of
$GL(m|n)$ in positive characteristic, extending work of Donkin \cite{D1}.
The answer is a natural generalization of the ``hook theorem''
of Berele and Regev \cite{BR1} and Sergeev \cite{S2} 
in characteristic $0$.

\vspace{2mm}
\noindent{\bf Ackowledgements.}
The idea that $GL(n|n)$ could be used to prove the 
Mullineux conjecture
was inspired by a paper of A. Regev \cite{R1}. We are especially grateful to
V. Serganova for explaining the material in $\S$\ref{S:serg}
to us during a visit
to the University of Oregon.
We would also like to thank A. Kleshchev for pointing out the reference
\cite{BOX}.

\section{The supergroup $GL(m|n)$}\label{S:basics}

Throughout, let $k$ be a field of 
characteristic $p \neq 2$.  All objects (superalgebras, supergroups, \dots) 
will be defined over $k$.  
A \emph{commutative superalgebra} is a $\mathbb{Z}_2$-graded associative
algebra $A = A_{\0} \oplus A_{\1}$
with $ab=(-1)^{\bar a \bar b}ba$ for all homogeneous $a, b \in A$,
where 
$\bar x \in \mathbb{Z}_2$ 
denotes the parity of a homogeneous vector $x$ in a vector superspace.
For an account of the basic language of superalgebras and
supergroups adopted here, we refer
the reader to \cite{B1, BK1}, see also \cite{J1}, \cite{K},
\cite[ch.I]{Leites} and \cite[ch.3, $\S\S$1--2, ch.4, $\S$1]{Man}.

The supergroup $G = GL(m|n)$ is the functor
from the category of commutative superalgebras to the category of groups
defined on a commutative superalgebra $A$ 
by letting $G(A)$ be the group of all  invertible $(m+n) \times (m+n)$
matrices of the form 
\begin{equation}\label{E:matrix}g=\left( 
\begin{array}{l|l}
W&X\\\hline
Y&Z
\end{array}
\right)
\end{equation}
where $W$ is an $m \times m$ matrix with entries in 
$A_{\0}$, $X$ is an $m\times n$ matrix with entries in 
$A_{\1}$, $Y$ is an $n\times m$ with entries in $A_{\1}$, and $Z$ is an $n \times n$ matrix with entries in $A_{\overline{0}}$.  If $f: A \to B$ is a superalgebra homomorphism, then $G(f):G(A) \to G(B)$ is the group homomorphism defined by applying $f$ to the matrix entries.  

Let $Mat$ be the affine superscheme with $Mat(A)$ consisting of {all}
(not necessarily invertible)
$(m+n)\times(m+n)$ matrices of the above form.
For $1 \leq i,j \leq m+n$, let $T_{i,j}$ be the function mapping a matrix
to its $ij$-entry. 
Then, the coordinate ring $k[Mat]$ is the free
commutative superalgebra on the generators
$\{T_{i,j}\:|\:1 \leq i, j \leq m+n\}$. Writing $\bar i = \0$ for
$i = 1,\dots,m$ and $\bar i = \1$ for $i = m+1,\dots,m+n$, the
parity of the generator $T_{i,j}$ is $\bar i + \bar j$.
By \cite[I.7.2]{Leites},
a matrix $g \in Mat(A)$ of the form (\ref{E:matrix}) is invertible if and 
only if $\det W \det Z \in A^\times$.
Hence, $G$ is the principal open subset of $Mat$ defined by the function
$\det:g \mapsto \det W \det Z$. In particular, the coordinate ring
$k[G]$ is the localization of $k[Mat]$ at $\det$.

Just like for
group schemes \cite[I.2.3]{J1}, the coordinate ring
$k[G]$ has the naturally induced structure of a 
Hopf superalgebra. Explicitly, the
comultiplication and counit of $k[G]$ are the unique 
superalgebra maps satisfying
\begin{align}\label{hstruct}
\Delta(T_{i,j}) &= \sum_{h =1}^{m+n} T_{i,h} \otimes T_{h,j},\phantom{(-1)^{(\bar i + \bar h)(\bar h + \bar j)}} \\
\varepsilon(T_{i,j})&=\delta_{i,j}\label{hhstruct}
\end{align}  
for all $1 \leq i,j \leq m+n$.
The subalgebra $k[Mat]$ of $k[G]$
is a subbialgebra but not a Hopf subalgebra, as it is not
invariant under the antipode.

It is sometimes convenient to work with an alternative set of generators
for the coordinate ring $k[G]$: define
\begin{equation}\label{altt}
\tilde T_{i, j} = (-1)^{\bar i (\bar i + \bar j)} T_{i,j}.
\end{equation}
In terms of these new generators, 
(\ref{hstruct}) becomes
\begin{align}\label{hstruct2}
\Delta(\tilde T_{i,j}) &= \sum_{h =1}^{m+n}(-1)^{(\bar i + \bar h)(\bar h + \bar j)} \tilde T_{i,h} \otimes \tilde T_{h,j}.
\end{align}  

A representation of $G$ means 
a natural transformation $\rho:G \rightarrow GL(M)$
for some vector superspace $M$, where $GL(M)$ is the supergroup
with $GL(M)(A)$ being equal to the group of all even automorphisms of 
the $A$-supermodule $M \otimes A$, for each commutative superalgebra $A$.
Equivalently, as with group schemes \cite[I.2.8]{J1}, $M$ is  
a right $k[G]$-comodule, i.e. there is an even
structure map $\eta: M \to M \otimes k[G]$ satisfying the usual comodule 
axioms. 
We will usually refer to such an $M$ as a {\em $G$-supermodule}.
For example, we have the \emph{natural representation} $V,$ 
the $m|n$-dimensional vector superspace with canonical basis
$v_1,\dots,v_m,v_{m+1},\dots,v_{m+n}$ where $\bar v_i = \bar i$.
Identify elements of $V \otimes A$ with column vectors via
\begin{equation*}
\sum_{i = 1}^{m+n} v_{i}\otimes a_{i} \longleftrightarrow 
\left(\begin{matrix} a_{1} \\
\vdots \\
a_{m+n}
\end{matrix} \right).
\end{equation*}
Then, the $G(A)$-action on $V \otimes A$ is the usual one by left multiplication. The induced comodule structure map $\eta:V \to V \otimes k[G]$ 
is given explicitly by 
\begin{equation}\label{vstruct}
\eta(v_{j}) = \sum_{i =1}^{m+n}v_{i} \otimes T_{i,j}
=
\sum_{i=1}^{m+n} (-1)^{\bar i (\bar i + \bar j)}
v_i \otimes \tilde T_{i,j}.
\end{equation}

The underlying purely even group $G_{\ev}$ of $G$ is by definition
the functor from superalgebras to groups 
with $G_\ev(A) := G(A_{\0})$. Thus, $G_\ev(A)$ consists
of all invertible matrices of the form (\ref{E:matrix}) with $X=Y=0$,
so $G_\ev \cong GL(m) \times GL(n)$.
Let $T$ be the usual maximal torus of $G_\ev$ consisting of diagonal matrices.
The character group $X(T) = \Hom(T, \mathbb{G}_m)$ is the free abelian
group on generators $\varepsilon_1,\dots,\varepsilon_m,\varepsilon_{m+1},\dots,\varepsilon_{m+n}$, where $\varepsilon_i$
picks out the $i$th diagonal entry of a diagonal matrix.
Put a symmetric bilinear form on $X(T)$ by declaring that
\begin{equation}\label{bf}
(\varepsilon_i, \varepsilon_j)
= (-1)^{\bar i} \delta_{i,j}.
\end{equation}
Let $W \cong S_m \times S_n$ 
be the Weyl group of $G_{\ev}$ with respect to $T$, identified
with the subgroup of $G_\ev$ consisting of 
all permutation matrices.

A \emph{full flag} $F = (F_1\subset \dots \subset F_{m+n})$ 
in the vector superspace $V$ means a 
chain of subsuperspaces 
of $V$
with each $F_{i}$ having dimension $i$ as a vector space.  
If $(u_{1}, u_2, \dotsc , u_{m+n})$ is an ordered homogeneous basis for $V$, 
we write $F(u_{1},u_2,\dotsc ,u_{m+n})$ for the full flag 
with $F_i = \langle u_1,\dots,u_i \rangle$.
By definition, a {\em Borel subgroup} $B$ of $G$ is the stabilizer of a full 
flag $F$ in $V$, i.e. $B(A)$ is the stabilizer in 
$G(A)$ of the canonical image of $F$ in $V \otimes A$ 
for each commutative superalgebra $A$.
Since $GL(m)$ (resp. $GL(n)$) acts transitively on the bases of 
$V_{\0}$ (resp. $V_{\1}$), it is easy to see two full flags 
$F$ and $F'$ in $V$ are conjugate under $G$ 
if and only if the superdimension of $F_i$ equals the superdimension of
$F_i'$ for each $i = 1,\dots,m+n$.
Consequently there are $\binom{m+n}{n}$ different
conjugacy classes of Borel subgroups.  

View the Weyl group $W$ of $G$
as the parabolic subgroup $S_m \times S_n$
of the symmetric group $S_{m+n}$ in the obvious way.
Let $D_{m,n}$ be the set of all minimal length
$S_m \times S_n \backslash  S_{m+n}$-coset representatives,
i.e. 
$$
D_{m,n} = \{w \in S_{m+n}\:|\:w^{-1} 1 < \dots < w^{-1} m,
w^{-1}(m+1) < \dots < w^{-1}(m+n)\}.
$$
For $w \in S_{m+n}$, let $B_w$ be the stabilizer of the full flag
$F(v_{w 1}, v_{w 2}, \dots, v_{w (m+n)})$.
Then, the Borel subgroups $\{B_w\:|\:w \in D_{m,n}\}$
give a set of representatives for the conjugacy classes of Borel 
subgroup in $G$ (cf. \cite[Proposition 1.2(a)]{Kacnote}). 
We point out that for $w \in D_{m,n}$, the underlying
even subgroup of $B_w$ is always the usual upper triangular
Borel subgroup $B_{\ev}$ of $G_{\ev}$.

The root system of $G$ is the set $\Phi = \{\varepsilon_i - \varepsilon_j\:|\:1 \leq i,j \leq m+n, i\neq j\}$. There are even and odd roots, the 
parity of the root 
$\varepsilon_i - \varepsilon_j$ being $\bar i + \bar j$.
Choosing $w \in S_{m+n}$ fixes a choice $B_w$ of 
Borel subgroup of $G$ containing $T$, hence a set
\begin{equation}\label{proots}
\Phi_w^+ = \{ \varepsilon_{w i} - \varepsilon_{w j}\:|\:1 \leq i < j \leq m+n\}
\end{equation}
of positive roots.
The corresponding dominance ordering on $\Phi$ is denoted $\leq_w$,
defined by
$\lambda \leq_w \mu$ if $\mu - \lambda \in \mathbb{Z}_{\geq 0} \Phi_w^+$.

For examples, first take $w =1$.
Then,
$B_1 = \operatorname{stab}_G F(v_1,v_2,\dots,v_{m+n})$ 
is the Borel subgroup with 
$B_1(A)$ consisting of all upper triangular invertible matrices
of the form (\ref{E:matrix}).
This is the {\em standard choice} of Borel subgroup,
giving rise to the standard choice of positive roots
$\Phi_1^+$ and the standard dominance ordering $\leq_1$ on $X(T)$.
Instead, let $w_0$ be the longest element of $S_m \times S_n$ and 
$w_1$ be the longest element of $D_{m,n}$, so that
$w_0 w_1$ is the longest element of the symmetric group $S_{m+n}$.
Then,
$B_{w_1} = \operatorname{stab}_G F(v_{m+1}, \dots, v_{m+n}, v_1, \dots, v_m)$ 
is the Borel
with $B_{w_1}(A)$ consisting of all invertible matrices of the 
form (\ref{E:matrix}) with $X=0$ and $W, Z$ upper triangular.
Finally, $B_{w_0 w_1} = w_0 {B_{w_1}} 
w_0^{-1}$ is the Borel subgroup of all lower triangular matrices.

\section{The superalgebra of distributions}\label{S:algebraofdistributions}

We next recall the definition 
of the superalgebra of distributions $\dist{}$ of $G$, 
following \cite[$\S$4]{B1}.
Let $\I$ be the kernel of the counit $\varepsilon:k[G] \rightarrow k$,
a superideal of $k[G]$. For $r \geq 0$, let 
\begin{align*}
\dist{r}&=\{x \in k[G]^{*}\:|\: x(\I^{r+1})=0 \} \cong (k[G]/\I^{r+1})^{*},\\
\dist{}&=\bigcup_{r \geq 0} \dist{r}.
\end{align*}
There is a multiplication on $k[G]^{*}$ dual to the comultiplication on $k[G]$, defined by
$(x y)(f) = (x \bar\otimes y) (\Delta(f))$ for $x, y \in k[G]^{*}$ and 
$f \in k[G]$. Note here (and later on)
we are implicitly using the superalgebra
rule of signs: $(x \bar \otimes y)(f \otimes g) =
(-1)^{\bar y \bar f} x(f) y(g)$.
One can check that $\dist{}$ is a subsuperalgebra of $k[G]^{*}$ 
using the fact that for $f \in \I$,
\begin{equation*}\label{E:lieeq}
\Delta(f) \in 1 \otimes f + f \otimes 1 + \I \otimes \I, 
\end{equation*}
or, more generally,
\begin{equation}\label{E:filteredeq}
\Delta(f_{1} \dotsb f_{r}) \in \prod_{i=1}^{r}(1 \otimes f_{i} + f_{i} \otimes 1) + \sum_{j=1}^{r} \I^{j}\otimes \I^{r+1-j}
\end{equation}
for all $f_{1},\dotsc ,f_{r} \in \I $.  In fact, since 
$\I^{r+1} \subseteq \I^{r}$, we have $\dist{r} \subseteq \dist{r+1}$ and 
(\ref{E:filteredeq}) shows that $\dist{r} \dist{s} 
\subseteq \dist{r+s}$, i.e. $\dist{}$ is a filtered superalgebra.  
By (\ref{E:filteredeq}) again, the subspace
$$
T_{1}(G) = \{x \in \dist{1}\:|\: x(1)=0 \} \cong (\I / \I^{2})^{*}
$$  
is closed under the superbracket
$[x,y] := xy - (-1)^{\bar x \bar y} yx$,
giving $T_1(G)$ the structure of Lie superalgebra, denoted
$\lie$. 
Finally, given a $G$-supermodule $M$ with structure map
$\eta:M \rightarrow M \otimes k[G]$, we can view $M$ as a $\dist{}$-supermodule
by $x.m = (1 \bar\otimes x) (\eta(m))$. In particular, this makes $M$ into
a $\lie$-supermodule. 

To describe $\lie$ explicitly in our case, recall
the alternative generators $\tilde T_{i,j}$ of $k[G]$ from
(\ref{altt}).
The superideal $\I$ is generated by 
$\{\tilde T_{i,j} - \delta_{i,j}\:|\:1 \leq i,j \leq m+n\}.$
So $\lie$ has a unique basis $\{e_{i,j}\:|\:1 \leq i,j \leq m+n\}$
such that $e_{i,j}(\tilde T_{h,l}) = \delta_{i,h} \delta_{j,l}$.
The parity of $e_{i,j}$ is $\bar i + \bar j$,
while (\ref{hstruct}) implies that the multiplication satisfies
\begin{equation}\label{meq}
[e_{i,j}, e_{h,l}] = \delta_{j,h} e_{i,l} - (-1)^{(\bar i + \bar j)(\bar h+\bar l)} \delta_{i,l} e_{h,j}.
\end{equation}
Thus $\lie$ is identified with the Lie superalgebra 
$\mathfrak{gl}(m|n)$ over $k$, see \cite{K}, 
so that $e_{i,j}$ corresponds to the $ij$-matrix unit.
By (\ref{vstruct}), the induced action of $\lie$ on the
natural representation $V$ of $G$ is given by 
$e_{i,j} v_h = \delta_{j,h} v_i$, i.e. $V$ 
is identified with the natural representation of $\mathfrak{gl}(m|n)$.

To describe $\dist{}$ explicitly, first note that over $\mathbb{C}$,
$\dist{}$ is simply the universal enveloping superalgebra of $\lie$.
To construct $\dist{}$ in general, 
let ${U}_{\mathbb{C}}$ be the universal enveloping superalgebra of 
the Lie superalgebra $\gl$ over $\mathbb{C}$.  
By the PBW theorem for Lie superalgebras (see \cite{K}), 
${U}_{\mathbb{C}}$ has basis consisting of all monomials
\[
\prod_{\substack{1 \leq i,j \leq m+n\\\bar i + \bar j = \0}}
e_{i,j}^{a_{i,j}}\prod_{\substack{1 \leq i,j \leq m+n\\\bar i+ \bar j=\1}}
e_{i,j}^{d_{i,j}}
\]
where $a_{i,j} \in \mathbb{Z}_{\geq 0}$, $d_{i,j} \in \{0,1 \}$, and the product is taken in any fixed order.  
We shall write $h_{i} =e_{i,i}$ for short.

Define the \emph{Kostant $\mathbb{Z}$-form} 
${U}_{\mathbb{Z}}$ to be the $\mathbb{Z}$-subalgebra of 
${U}_{\mathbb{C}}$ generated by elements
$e_{i,j}\:(1 \leq i,j \leq m+n, \bar i + \bar j = \1)$,
$e_{i,j}^{(r)}\: 
(1 \leq i \neq j \leq m+n, \bar i + \bar j = \0, r \geq 1)$
and
$ \binom{h_i}{r}\:(1 \leq i \leq m+n, r \geq 1)$.
Here, $e_{i,j}^{(r)} := e_{i,j}^r / (r!)$ and 
$\binom{h_{i}}{r} := {h_{i}(h_{i}-1) \dotsb (h_{i}-r+1)}/{(r!)}$.
Following the proof of \cite[Th.2]{S1}, one verifies the following:

\begin{lemma}\label{stbas}
The superalgebra ${U}_{\mathbb{Z}}$ is a $\mathbb{Z}$-free
$\mathbb{Z}$-module with basis being given by the set of all monomials
of the form
$$
\prod_{\substack{ 1 \leq i,j \leq m+n \\
                   \bar{i}+\bar{j} = \0}} e_{i,j}^{(a_{i,j})}
\prod_{1 \leq i \leq m+n} \binom{h_{i}}{r_{i}}
\prod_{\substack{1 \leq i, j \leq m+n \\
\bar{i}+\bar{j} = \1}} e_{i,j}^{d_{i,j}}
$$
for all $a_{i,j}, r_{i} \in \mathbb{Z}_{\geq 0}$ and $d_{i,j} 
\in \{0,1 \}$, where the product is taken in any fixed order.
\end{lemma}

The enveloping superalgebra ${U}_{\mathbb{C}}$ is a 
Hopf superalgebra in a canonical way, hence
${U}_{\mathbb{Z}}$ is a Hopf superalgebra over $\mathbb{Z}$.
Finally, set
${U}_k = k \otimes_{\mathbb{Z}} {U}_{\mathbb{Z}},$
naturally a Hopf superalgebra over $k$.
We will abuse notation by using the same symbols
$e_{i,j}^{(r)}, \binom{h_i}{r}$ etc... for the
canonical images of these elements of ${U}_{\mathbb{Z}}$ in
${U}_k$.
Now the basic fact is the following:

\begin{theorem}\label{T:Distiso}  ${U}_{k}$ and $\dist{}$ are isomorphic as Hopf superalgebras.
\end{theorem}

\begin{proof}
In the case when $k = \mathbb{C}$, 
the isomorphism $i:{U}_{\mathbb{C}} \rightarrow \dist{}$
is induced by the Lie superalgebra isomorphism
mapping the matrix unit $e_{i,j} \in \gl$ to the element with the same
name in $\lie$.
For arbitrary $k$, the isomorphism
$i:{U}_k \rightarrow \dist{}$ is obtained by reducing
this one modulo $p$.
\end{proof}

In view of the theorem, we will henceforth 
\emph{identify} ${U}_{k}$ with $\dist{}$.
It is also easy to describe the superalgebras of distributions of
our various natural subgroups of $G$ as subalgebras
of $\dist{}$.
For example, $\Dist{T}$ is the subalgebra generated by all
$ \binom{h_i}{r}\:(1 \leq i \leq m+n, r \geq 1)$,
$\Dist{B_{\ev}}$ is the subalgebra generated by $\Dist{T}$
and all
$e_{i,j}^{(r)}\: 
(1 \leq i < j \leq m+n, \bar i + \bar j = \0, r \geq 1)$,
and for $w \in D_{m,n}$,
$\Dist{B_w}$ is the subalgebra generated by $\Dist{B_{\ev}}$
and all $e_{i, j} \:(1 \leq i, j \leq m+n, \bar i + \bar j = 1,
w^{-1} i < w^{-1} j).$

For $\lambda=\sum_{i=1}^{m+n} \lambda_{i}\varepsilon_{i} \in X(T)$  
and a $\dist{}$-supermodule $M$, define the {\em $\lambda$-weight space of $M$}
to be
\begin{equation}\label{wtsp}
M_{\lambda}=\left\{m \in M\:\bigg|\:
\binom{h_{i}}{r}m=\binom{\lambda_{i}}{r}m \mbox{ for all } i =1,\dots,m+n,
r \geq 1 \right\}. 
\end{equation}
We call a $\dist{}$-supermodule $M$ {\em integrable}
if it is locally finite over $\dist{}$ and satisfies
$M = \sum_{\lambda \in X(T)} M_{\lambda}$.
If $M$ is a $G$-supermodule viewed as a $\dist{}$-supermodule in 
the natural way, then $M$ is integrable.
The goal in the remainder of the section is to prove 
conversely that any integrable $\dist{}$-supermodule
can be lifted in a unique way to $G$.

Let $\dist{}^\diamond$ denote the
\emph{restricted dual} of $\dist{}$, namely,
the set of all $f \in \dist{}^{*}$ such that $f(I) = 0$ 
for some two-sided superideal $I \subset \dist{}$ (depending on $f$) with
$\dist{} /I$ being a finite dimensional integrable $\dist{}$-supermodule.  
If $M$ is an integrable $\dist{}$-supermodule with homogenous basis 
$\{m_{i} \}_{i\in I}$, 
its 
\emph{coefficent space} $cf(M)$
is the subspace of $\dist{}^{*}$ spanned by the 
\emph{coefficient functions} $f_{i,j}$ defined by
\begin{equation}
\label{cfd}
um_{j} = 
(-1)^{\bar u \bar m_j} \sum_{i \in I} f_{i,j}(u) m_i
\end{equation}
for all homogeneous $u \in \dist{}$.  
Note that this definition is independent of the choice of homogenous basis.  
As in the purely even case \cite[(3.1a)]{Donk}, 
we have the following lemma:

\begin{lemma}\label{L:restricteddual}  $f \in \dist{}^{*}$ belongs to $\dist{}^{\diamond}$ if and only if $f \in cf(M)$ for some integrable $\dist{}$-supermodule $M$.
\end{lemma}

If $M$ and $N$ are integrable $\dist{}$-supermodules, then 
$M \otimes N$ is also an integrable supermodule and $cf(M \otimes N) = cf(M)cf(N)$.  Consequently, Lemma \ref{L:restricteddual} implies 
$\dist{}^{\diamond}$ is a subsuperalgebra of $\dist{}^{*}$.  
Indeed, $\dist{}^{\diamond}$ has a natural Hopf superalgebra structure
dual to that on $\dist{}$, cf. the argument after \cite[Lemma 5.2]{B1}.

\begin{theorem}\label{it}
The map $\iota:k[G] \rightarrow \dist{}^\diamond$ defined
by $\iota(f)(u)=(-1)^{\bar{f}\bar{u}}u(f)$ for all 
homogeneous $f \in k[G]$ and $u \in \dist{}$
is an isomorphism of Hopf 
superalgebras.
\end{theorem}

\begin{proof}
Note $\iota$ is automatically a Hopf superalgebra homomorphism,
since the Hopf superalgebra structure on $\dist{}$ is 
dual to that on $k[G]$ and the Hopf superalgebra structure on
$\dist{}^\diamond$ is dual to that on $\dist{}$.
Furthermore if $\iota({f})=0$ then $u(f)=0$ for all $u \in \dist{r}$, 
so $f \in \I^{r+1}$.  
Since $r$ was arbitrary we deduce $f \in \bigcap_{r \geq 0} \I^{r+1}$,
hence $f = 0$. This shows that $\iota$ is injective.
It remains to prove that $\iota$ is surjective.

Fix an order for the products in the monomials in the PBW basis for 
$\dist{}$ from Lemma~\ref{stbas} 
so that all monomials are of the form $mu$ where $m$ is a monomial in the 
$e_{i,j}$ with $\bar{i} +\bar{j}=1$ and $u \in \distev{}$.  
Let $\Gamma = \{(i,j): 1 \leq i,j \leq m+n, \bar{i}+\bar{j}=\bar{1} \}$.  
For each $I \subseteq \Gamma$, let $m_{I}$ denote the 
PBW monomial given by taking the product of the 
$e_{i,j}$'s for $(i,j) \in I$ in the fixed order.   
By Lemma \ref{stbas} we have the vector space decomposition 
\[
\dist{}= \bigoplus_{I \subseteq \Gamma} m_{I}\distev{}.
\] 
For $I \subseteq \Gamma$, let $\eta_{I} \in \dist{}^{*}$ be the 
linear functional given by $\eta_{I}(m_{I})=1$ and $\eta_{I}(m)=0$ 
for any other ordered PBW monomial different from $m_{I}$.

\vspace{1mm}
\noindent{\bf Claim 1.} {\em
For any $I \subseteq \Gamma$, we have that $\eta_{I} \in \iota(k[G])
\subseteq \dist{}^\diamond$.}
\vspace{1mm}

To prove this, let $N=m^{2}+n^{2}$.  Let 
$M$ denote $\bigwedge^{N}(V \otimes V^{*})$
viewed as a $\dist{}$-supermodule in the natural way.
Since $M$ is in fact a $G$-supermodule, we 
have that $cf(M) \subseteq \iota(k[G])$.  
Therefore to prove Claim 1, it suffices to 
show that $\eta_{I} \in cf(M)$ for any $I \subseteq \Gamma$.
Let $f_1,\dots,f_{m+n}$ be the basis for $V^*$ dual to the basis
$v_1,\dots,v_{m+n}$ of $V$.
Let $z_{i,j} = v_{i} \otimes f_{j} \in V \otimes V^{*}$.  
Fix a total order on the set $\{1,\dots,m+n\}\times\{1,\dots,m+n\}$ 
and in this order let 
$\Sigma$ be the set of all weakly increasing sequences 
$S=((i_{1},j_{1}) \leq \dotsb \leq (i_{N},j_{N}))$ of length $N$ such that 
$(i_{k},j_{k}) < (i_{k+1},j_{k+1})$ whenever $\bar{i}_k+\bar{j}_k=\0$.  
For $S \in \Sigma$, let $z_{S} = 
z_{i_{1},j_{1}} \wedge \dotsb \wedge z_{i_{N},j_{N}}$, so that
$\{z_{S}\}_{S \in \Sigma}$ is a basis for $M$.  In particular, 
let $z=z_{S}$ for the sequence $S$ containing all $(i,j)$ with 
$\bar{i}+\bar{j}=\0$.  
Then $z$ spans $\bigwedge^{N}((V \otimes V^{*})_{\bar{0}})=\bigwedge^{N}
(V_{\bar{0}} \otimes V^{*}_{\bar{0}} \oplus V_{\bar{1}} \otimes V^{*}_{\bar{1}})$, which is a 1-dimensional trivial $\distev{}$-submodule of $M$.

Observe now that 
$\{m_{I}z\}_{I \subseteq \Gamma}$ is a linearly independent set of 
homogeneous vectors, because they are related to the basis elements
$\{z_{S}\}_{S \in \Sigma}$ in a unitriangular way.
Extend this set to a homogeneous basis $\mathscr{B}$ of $M$.  
For $I \subseteq \Gamma$ and $u \in \dist{}$ define $g_{I}(u)$ to be the 
$m_{I}z$ coefficent of $uz$ when expressed in the basis $\mathscr{B}$.  
Then $g_{I}(m_{J})=\delta_{I,J}$ for all $I, J \subseteq \Gamma$.  
Furthermore, since $z$ spans a trivial $\distev{}$-module, $uz=0$ for all 
monomials in our ordered PBW basis for $\dist{}$ not of the form $m_{J}$, i.e.
$g_{I}(u)=0$ for all such monomials.  Therefore 
$\eta_{I}=g_{I} \in cf(M)$, proving the claim.

\vspace{1mm}
\noindent{\bf Claim 2.} {\em
For any $I \subseteq \Gamma$ and $f \in \distev{}^\diamond$,
we have that
$(\eta_I f)(m_I u) = f(u)$ and $(\eta_I f)(m_J u) = 0$
for all $u \in \distev{}$ and $J \nsupseteq I$.
}
\vspace{1mm}

Indeed, by the definition of multiplication in $\dist{}^{\diamond}$, 
we have $(\eta_{I}f)(m_{J}u)=(\eta_{I} \bar{\otimes} f)
(\delta(m_{J}u))$, where $\delta$ is the comultiplication on $\dist{}$.  
Recalling that $\delta(e_{i,j}) = e_{i,j} \otimes 1 + 1 \otimes e_{i,j}$,
we see that, when expressed in the ordered PBW basis of 
$\dist{} \otimes \dist{}$, the $(m_{I} \otimes -)$-component of 
$\delta(m_{J}u)$ is equal to $m_{I}\otimes u$ if $J=I$ and $0$ if 
$J \nsupseteq I$.  This implies the claim.

\vspace{1mm}
\noindent{\bf Claim 3.} {\em
For any $f \in \dist{}^\diamond$ and $I \subseteq \Gamma$, 
there is a function
$f_I \in \iota(k[G])$ such that $f_I = f$ on $m_I \distev{}$ and
$f_I = 0$ on $\bigoplus_{J \not\supseteq I} m_J \distev{}$.
}
\vspace{1mm}

To prove this, we need 
to appeal to the analogous theorem for the underlying
even group $G_{\ev}$.
Just as for $\dist{}$ we can define integrable $\Dist{G_{\ev}}$-supermodules, 
coefficent space, the restricted dual $\distev{}^\diamond$, etc...
By the purely even theory, 
the natural map $\iota_{\ev}:
k[G_{ev}] \rightarrow \distev{}^\diamond$ (the analogue of the map
$\iota:k[G] \rightarrow \dist{}^\diamond$ being considered here)
is an isomorphism,
see e.g. \cite[(3.1c)]{Donk} for the proof.
An integrable $\dist{}$-supermodule 
is integrable over $\distev{}$ too, so restriction gives us a Hopf 
superalgebra homomorphism
$\vartheta: \dist{}^{\diamond} \to \distev{}^{\diamond}$  
such that
$\vartheta \circ \iota = \iota_{\ev} \circ \varphi$, where
$\varphi:k[G] \twoheadrightarrow k[G_\ev]$ is the canonical map
induced by the inclusion of $G_{\ev}$ into $G$.

Now take $f \in \dist{}^\diamond$ and write 
$\delta(f) = \sum_{j} f_{j}\otimes g_{j}$. By the previous paragraph,
we can find even elements $h_j \in \iota(k[G])$ such that 
$\vartheta(g_{j}) = \vartheta(h_{j})$ for each $j$.  
For $I \subseteq \Gamma$, 
let $f_{I}=\sum_{j} f_{j}(m_{I})\eta_{I}h_{j}$, an element of
$\iota(k[G])$ by Claim 1.  
By Claim 2, we have $f_{I}=f$ on $m_{I}\distev{}$ and $f_{I}=0$ on $\sum_{J \nsupseteq I} m_{J}\distev{}$, as required to prove the claim.

Now we can complete the proof.
Fix $f \in \dist{}^{\diamond}$.  
For $i=0,1, \dotsc , 2mn$ define $f^{(i)}$ recursively by 
\begin{align*}
f^{(0)}=f-f_{\varnothing} & & f^{(i)}=f^{(i-1)}-\sum_{I \subseteq \Gamma, |I|=i} (f^{(i-1)})_{I},
\end{align*}
invoking Claim 3.
An easy induction on $i$ using Claim 3 shows that 
$f^{(i)}=0$ on $\bigoplus_{J \subseteq \Gamma, |J| \leq i} m_{J}\distev{}$.  
In particular, $f^{(2mn)}=0$ on $\dist{}$.  This implies the surjectivity
of $\iota$, since $f$ is obtained from $f^{(2mn)}$ by adding elements
of $\iota(k[G])$.
\end{proof}

\begin{corollary}\label{C:vermaconjecture}  
The category of $G$-supermodules is isomorphic to the category of integrable $\dist{}$-supermodules.
\end{corollary}

\begin{proof}  
Say $M$ is an integrable $\dist{}$-supermodule with homogenous basis 
$\{m_{i}\}_{i \in I}$.  
Let $f_{i,j}$ be the corresponding coefficent functions
defined according to (\ref{cfd}).
By Theorem~\ref{it}, there are unique $g_{i,j} \in k[G]$ such that 
$\iota(g_{i,j}) = f_{i,j}$. Define a structure map
$\eta: M \to M \otimes k[G]$
making $M$ into a $G$-supermodule by 
\[
\eta(m_j) = \sum_{i \in I} m_{i} \otimes g_{i,j}.
\]
Conversely, as discussed at the beginning of the section,
any $G$-supermodule can be viewed as an integrable $\dist{}$-supermodule in
a natural way.  One 
can verify that these two constructions give mutually inverse
functors between the two categories.
\end{proof}

In view of the corollary, we will not distinguish between
$G$-supermodules and integrable $\dist{}$-supermodules in the rest
of the article.

\section{Highest weight theory}\label{S:serg}

Now we describe the classification of the irreducible representations
of $G$ by their highest weights. It seems to be more convenient to
work first in the category $\mathscr O$ of all $\dist{}$-supermodules $M$ 
that are locally finite over $\Dist{B_{\ev}}$ and
satisfy $M = \bigoplus_{\lambda \in X(T)} M_\lambda$.
Fix a choice of $w \in D_{m,n}$, hence a Borel subgroup $B_w$
and dominance ordering $\leq_w$ on $X(T)$.
By Lemma~\ref{stbas}, $\Dist{B_w}$ is a
free right $\Dist{B_{\ev}}$-module of finite rank.
So the condition 
that $M$ is locally finite over
$\Dist{B_{\ev}}$ 
in the definition of category $\mathscr O$
is equivalent to $M$ being locally finite
over $\Dist{B_w}$.
For $\lambda \in X(T)$, we have the {\em Verma module}
$$
M_w(\lambda) := \Dist{G} \otimes_{\Dist{B_w}} k_\lambda,
$$
where $k_\lambda$ denotes $k$ viewed as a 
$\Dist{B_w}$-supermodule of weight $\lambda$.
We say that a vector $v$ in a $\Dist{G}$-supermodule $M$
is a {\em $w$-primitive vector of weight $\lambda$} if
$\Dist{B_w} v \cong k_\lambda$ as a $\Dist{B_w}$-supermodule.
Familiar arguments exactly as for semisimple Lie algebras
over $\mathbb C$ show:

\begin{lemma}\label{hw} Let $w \in D_{m,n}$ and $\lambda \in X(T)$.
\begin{itemize}
\item[(i)]
The $\lambda$-weight space of $M_w(\lambda)$ is $1$-dimensional,
and all other weights of $M_w(\lambda)$ are $<_w \lambda$.
\item[(ii)]
Any non-zero quotient of $M_w(\lambda)$
is generated by a $w$-primitive vector of weight $\lambda$,
unique up to scalars.
\item[(iii)]
Any $\dist{}$-supermodule generated by a $w$-primitive
vector of weight $\lambda$ is isomorphic to a quotient of 
$M_w(\lambda)$.
\item[(iv)]
$M_w(\lambda)$ has a unique irreducible quotient $L_w(\lambda)$,
and the $\{L_w(\lambda)\}_{\lambda \in X(T)}$ give
a complete set of pairwise non-isomorphic irreducibles in
$\mathscr O$.
\end{itemize}
\end{lemma}

In this way, we get a parametrization of the irreducible objects
in $\mathscr O$ by their highest weights with respect to the ordering
$\leq_w$.
Of course, the parametrization depends on the initial choice of 
$w \in D_{m,n}$.
To translate between labelings
arising from different choices $w, w' \in D_{m,n}$,
it suffices to consider the situation that $w, w'$ are adjacent
with respect to the usual Bruhat ordering on $D_{m,n}$.
In that case the following theorem of Serganova \cite{sergthesis},
see also \cite[Lemma 0.3]{PSI},
does the job. For the statement, recall the definition of the
form $(.,.)$ on $X(T)$ from (\ref{bf}).

\begin{lemma}\label{st}
Let $\lambda \in X(T)$.
Suppose that $w, w' \in D_{m,n}$ are 
adjacent in the Bruhat ordering, so 
$\Phi_{w'}^+ = \Phi_w^+ - \{\alpha\} \cup \{-\alpha\}$
for some odd root $\alpha = \varepsilon_i - \varepsilon_j \in \Phi$.
Then, 
$$
L_w(\lambda) \cong
\left\{
\begin{array}{ll}
L_{w'}(\lambda)&
\hbox{if $(\lambda,\alpha) \equiv 0 \pmod{p}$,}\\
L_{w'}(\lambda - \alpha)&\hbox{if $(\lambda, \alpha) \not\equiv 0 \pmod{p}$.}
\end{array}
\right.
$$
\end{lemma}

\begin{proof}
Let $v$ be a $w$-primitive vector in $L_w(\lambda)$ of weight $\lambda$,
cf. Lemma~\ref{hw}.
We claim first that $e_{r,s} e_{j,i} v = 0$ for all $1 \leq r,s \leq m+n$
with $\varepsilon_r - \varepsilon_s \in \Phi_w^+ \cap \Phi_{w'}^+$.
We know that $e_{r,s} v = 0$ as $v$ is
$w$-primitive. So we are done immediately if $[e_{r,s},e_{j,i}] = 0$.
In view of (\ref{meq}), this just leaves the possibilities $s = j$ or $r = i$.
Suppose first that $s = j$.
Noting that $w' = (i\:\:j) w$, the assumption that
$\varepsilon_r - \varepsilon_j \in \Phi_{w'}^+$ implies by (\ref{proots}) 
that $\varepsilon_r - \varepsilon_i \in \Phi_{w}^+$,
hence $e_{r,i} v = 0$. Therefore
$e_{r,j} e_{j,i} v = e_{r,i} v = 0$.
The remaining case when $r = i$ is similar.

Now suppose that $e_{j,i} v \neq 0$.
Since $e_{j,i}^2 = 0$, we get from the previous paragraph
that $e_{j,i} v$ is $w'$-primitive of weight $\lambda - \alpha$.
Hence, $L_w(\lambda) \cong L_{w'}(\lambda - \alpha)$.
On the other hand, if $e_{j,i} v = 0$, then $v$ itself is
already $w'$-primitive of weight $\lambda$ so
$L_w(\lambda) \cong L_{w'}(\lambda)$.
Thus, to complete the proof of the lemma, it suffices to show that
$e_{j,i} v \neq 0$ if and only if 
$(\lambda, \alpha) \not\equiv 0 \pmod{p}$. 
But $e_{j,i} v \neq 0$ if and only if 
there is some element $x\in\Dist{B_w}$ such that 
$x e_{j,i} v$ is a non-zero multiple of $v$.
In view of the first paragraph, the only $x$ that needs to be considered
is $e_{i,j}$. Finally, $e_{i,j} e_{j,i} v = (\lambda,\alpha) v$.
\end{proof}

Recall that $w_1$ denotes the longest element of $D_{m,n}$.
For $\lambda \in X(T)$, define $\widetilde{\lambda} \in X(T)$
from the isomorphism
\begin{equation}\label{Tdef}
L_1(\lambda) \cong L_{w_1}(\widetilde{\lambda}).
\end{equation}
Lemma~\ref{st} implies the following algorithm for computing
$\widetilde{\lambda}$:

\begin{theorem}\label{salg}
Pick an ordering $\beta_1, \dots, \beta_{mn}$
of the roots $\{\varepsilon_i - \varepsilon_j\:|\:
1 \leq i \leq m, m+1 \leq j \leq m+n\}$
such that $\beta_i \leq_1 \beta_j$ implies $i \leq j$.
Set $\lambda^{(0)} = \lambda$, and inductively define
$$
\lambda^{(i)} = \left\{
\begin{array}{ll}
\lambda^{(i-1)}&\hbox{if $(\lambda^{(i-1)}, \beta_{i}) \equiv 0 \pmod{p}$,}\\
\lambda^{(i-1)} - \beta_i&\hbox{if $(\lambda^{(i-1)}, \beta_{i}) 
\not\equiv 0 \pmod{p}$,}\\
\end{array}
\right.
$$
for $i = 1,\dots,mn$.
Then, $\widetilde{\lambda} = \lambda^{(mn)}$.
\end{theorem}

We refer to the algorithm for $\widetilde{\lambda}$
given by the theorem as {\em Serganova's algorithm}.
For an example, suppose $m = n = 2, p = 3$ and $\lambda = 
\varepsilon_1 + \varepsilon_2 + 2\varepsilon_3$.
Taking $\beta_1=
\varepsilon_2 - \varepsilon_3,
\beta_2=
\varepsilon_2 - \varepsilon_4,
\beta_3=
\varepsilon_1 - \varepsilon_3,\beta_4 = 
\varepsilon_1 - \varepsilon_4$, we get successively
$\lambda^{(1)} = \varepsilon_1 + \varepsilon_2 + 2 \varepsilon_3,
\lambda^{(2)} = \varepsilon_1 + 2\varepsilon_3 + \varepsilon_4,
\lambda^{(3)} = \varepsilon_1 + 2\varepsilon_3 + \varepsilon_4,
\lambda^{(4)} = 2\varepsilon_3 + 2\varepsilon_4$.
Hence, $\widetilde \lambda = 2 \varepsilon_3 + 2 \varepsilon_4$.

Now we pass from $\mathscr O$ to the finite dimensional irreducible
representations of $G$.
We will work now just with the standard choice of Borel subgroup
$B_1$ and the corresponding standard dominance ordering
$\leq_1$ on $X(T)$.
Let
$$
X^+(T) = \{\lambda = \sum_{i=1}^{m+n} \lambda_i \varepsilon_i \in X(T)\:|\:
\lambda_1 \geq \dots \geq \lambda_m, \lambda_{m+1} \geq \dots \geq \lambda_{m+n}\}
$$
denote the set of all {\em dominant integral weights}.
The proof of the first part of the following lemma goes back to 
Kac \cite{Kacnote}, while the second part is due to Serganova.

\begin{lemma}\label{fdclass}
Given any $\lambda \in X(T)$, $L_1(\lambda)$ is finite dimensional
if and only if $\lambda \in X^+(T)$.
Moreover, for $\lambda \in X^+(T)$, the $\leq_1$-lowest weight of
$L_1(\lambda)$ is $w_0 \widetilde \lambda$.
\end{lemma}

\begin{proof}
Suppose first that $L_1(\lambda)$ is finite dimensional for $\lambda \in X(T)$.
Then, it contains a $\Dist{B_{\ev}}$-primitive vector of weight
$\lambda$, hence by the purely even theory we must have that
$\lambda \in X^+(T)$. Conversely, suppose that
$\lambda \in X^+(T)$. Then, there is a finite dimensional irreducible
$\Dist{G_{\ev}}$-supermodule $L_{\ev}(\lambda)$ of highest weight $\lambda$.
Let $P$ be the closed subgroup of $G$ with $P(A)$ consisting of
all invertible matrices of the form (\ref{E:matrix}) with $Y = 0$.
We can view $L_{\ev}(\lambda)$ as a $\Dist{P}$-supermodule
so that all $e_{i,j}$ for $1 \leq i \leq m, m+1 \leq j \leq m+n$ act
as zero.
Consider the induced supermodule
$$
\Dist{G} \otimes_{\Dist{P}} L_{\ev}(\lambda).
$$
It is a finite dimensional module generated by a $1$-primitive vector
of weight $\lambda$. Hence, $M_1(\lambda)$ has a finite dimensional
quotient. This shows that $L_1(\lambda)$ is finite dimensional.
Finally, by (\ref{Tdef}), $L_1(\lambda) \cong L_{w_1}(\widetilde \lambda)$.
Hence, all its weights are $\leq_{w_1} \widetilde \lambda$.
Since $L_1(\lambda)$ is finite dimensional, 
the Weyl group $W$ acts by permuting weights. 
Hence we can act with $w_0$ to get that all its weights
are $\geq_1 w_0 \widetilde \lambda$.
\end{proof}

Lemmas~\ref{hw} and \ref{fdclass} show that 
$\{L_1(\lambda)\}_{\lambda \in X^+(T)}$ is
a complete set of pairwise non-isomorphic 
irreducible integrable $\dist{}$-supermodules.
In view of Corollary~\ref{C:vermaconjecture}, we can lift 
the $\dist{}$-supermodule $L_1(\lambda)$ for $\lambda \in X^+(T)$
uniquely to $G$.
We will denote the resulting irreducible $G$-supermodule
simply by $L(\lambda)$ from now on.
To summarize, using the second part of Lemma~\ref{fdclass} 
for the statement about $L(\lambda)^*$, we have shown:

\begin{theorem}\label{hwt}
The supermodules $\{L(\lambda)\}_{\lambda \in X^+(T)}$ form a complete
set of pairwise non-isomorphic irreducible $G$-supermodules.
Moreover, for $\lambda \in X^+(T)$, 
$L(\lambda)^* \cong L(-w_0 \widetilde \lambda)$.
\end{theorem}

\begin{remark}\label{rks}\rm
(i) The second part of 
Theorem~\ref{hwt} implies that the restriction of the map
$\sim$ from Theorem \ref{salg} gives a bijection
$\sim:X^+(T) \rightarrow X^+(T)$.

(ii) 
A weight $\lambda = \sum_{i=1}^{m+n} \lambda_i \varepsilon_i
\in X^+(T)$ is called {\em restricted} if either $p = 0$ or $p > 0$
and
$\lambda_i - \lambda_{i+1} < p$ for each $i = 1,\dots,m-1,m+1,\dots,m+n-1$. 
Assuming now that $p > 0$, 
let $X^+(T)_{\operatorname{res}}$ 
denote the set of all restricted $\lambda \in X^+(T)$.
Let $\Fr:G \rightarrow G_{\ev}$ be the Frobenius
morphism defined on $g \in G(A)$ by raising all
the matrix entries of $g$ to the power $p$, for each commutative superalgebra
$A$. Let $G_1 = \ker \Fr$ be the Frobenius kernel.
By a similar argument to 
\cite[6.4]{Borel}, the restriction of 
$L(\lambda)$ to $G_1$ remains irreducible
for all $\lambda \in X^+(T)_{\operatorname{res}}$, see \cite{thesis}.

(iii) Again for $p > 0$, 
there is an analogue for $G$ of the Steinberg tensor product theorem.
Given (ii), the proof is essentially the same as in \cite[$\S$9]{B1},
see \cite{thesis} for the details.
To state the result, let $L_{\ev}(\lambda)$
denote the irreducible $G_{\ev}$-supermodule
of highest weight $\lambda \in X^+(T)$ as in the proof of Lemma~\ref{fdclass}.
Inflating through the Frobenius morphism $\Fr:G \rightarrow G_{\ev}$,
we obtain an irreducible $G$-supermodule $\Fr^* L_{\ev}(\lambda)\cong
L(p \lambda)$.
In general, for $\lambda \in X^+(T)$, we can write
$\lambda = \mu + p \nu$ where $\mu \in X^+(T)_{\operatorname{res}}$
and $\nu \in X^+(T)$. Steinberg's tensor product theorem shows
that
\begin{equation}
L(\lambda) \cong L(\mu) \otimes \Fr^* L_{\ev}(\nu).
\end{equation}

(iv) Note for any $\lambda \in X^+(T)$, 
$\Fr^* L_{\ev}(\lambda)$ is trivial over $G_1$.
So (ii), (iii) show in particular that $L(\lambda)$ is irreducible
over $G_1$ if and only if $\lambda \in X^+(T)_{\operatorname{res}}$.
Given this, the second part of Theorem~\ref{hwt} implies that
the set $X^+(T)_{\operatorname{res}}$ is stable under the map $\sim$.
Finally, take $\lambda = \mu + p \nu$
where $\mu \in X^+(T)_{\operatorname{res}}$ and $\nu \in X^+(T)$, as in (iii).
Then,
\begin{align*}
L(-w_0 \widetilde{\lambda}) &\cong
L(\lambda)^* \cong L(\mu)^* \otimes \Fr^* (L_{\ev}(\nu)^*)\\
&\cong L(-w_0 \widetilde{\mu}) \otimes \Fr^* L_{\ev}(-w_0 \nu)
\cong L(-w_0(\widetilde \mu + p\nu)).
\end{align*}
Hence, 
$\widetilde \lambda = \widetilde{\mu} + p\nu$.
This reduces the problem of computing $\widetilde \lambda$
to the special case that $\lambda$ is restricted.
\end{remark}

\section{Polynomial representations}\label{S:schuralgebra}

In this section, we discuss polynomial representations of $G$
in the spirit of Green's monograph \cite{G1}.
Let $A(m|n)$ denote the subbialgebra $k[Mat]$ of $k[G]$,
so $A(m|n)$ is the free commutative superalgebra on
the generators $\{\tilde T_{i,j}\}_{1\leq i,j \leq m+n}$ 
from (\ref{altt}).
Obviously, $A(m|n)$ is $\mathbb Z$-graded by degree,
\begin{equation}\label{ad}
A(m|n) = \bigoplus_{d \geq 0} A(m|n,d).
\end{equation}
The subspace $A(m|n,d)$ is a finite dimensional 
subcoalgebra of $A(m|n)$.
A representation $M$ of $G$ is called a \emph{polynomial representation} 
(resp. a {\em polynomial representation of degree $d$})
if the comodule structure map $\eta:M \rightarrow M \otimes k[G]$
has image contained in $M \otimes A(m|n)$ (resp. in $M \otimes A(m|n,d)$).
For example, the $d$th tensor power
$V^{\otimes d}$ of the natural representation of $G$
is polynomial of degree $d$.
In general, a $G$-supermodule $M$ is polynomial of degree $d$ if it is
isomorphic to a direct sum of subquotients of $V^{\otimes d}$.

By \cite[Lemma 5.1]{BK1},
the decomposition (\ref{ad}) induces a decomposition of any 
polynomial representation into a direct sum of homogeneous polynomial
representations.
Moreover, the category of polynomial representations of degree $d$ is
isomorphic to the category of supermodules over the
{\em Schur superalgebra}
\begin{equation}\label{sdef}
S(m|n,d) := A(m|n,d)^*,
\end{equation}
where the superalgebra structure on $S(m|n,d)$ 
is the one dual to the coalgebra structure on $A(m|n,d)$.
Thus, the polynomial representation theory of $G$ reduces to studying
representations of the finite dimensional superalgebras $S(m|n,d)$
for all $d \geq 0$. The latter has been investigated recently 
over a field of positive characteristic by
Donkin \cite{D1}, see also \cite{Muir}.

Let $I(m|n, d)$ denote the set of all functions
from $\{1,\dots,d\}$ to $\{1,\dots,m+n\}$.
We usually view $\ui \in I(m|n,d)$ as a $d$-tuple
$(i_1,\dots,i_d)$ with entries in $\{1,\dots,m+n\}$.
In order to write down the various signs that will arise,
introduce the notation
$\epsilon_{\ui} = (\bar i_1, \dots,\bar i_d) \in \mathbb Z_2^d$,
for any $\ui \in I(m|n,d)$.
For tuples $\epsilon = (\epsilon_1,\dots,\epsilon_d), 
\delta = (\delta_1,\dots,\delta_d) \in \mathbb Z_2^d$ and $w \in S_d$, let
\begin{align}
\alpha(\epsilon, \delta) &= \prod_{1 \leq s < t \leq d}
(-1)^{\delta_s \epsilon_t},\\
\gamma(\epsilon, w) &=  \prod_{\substack{1 \leq s < t \leq d \\ w^{-1} s > w^{-1} t}} (-1)^{\epsilon_s \epsilon_t}.
\end{align}

The symmetric group $S_d$ acts on the right on $I(m|n,d)$
by composition of functions, i.e. $(i_1,\dots,i_d) \cdot w 
= (i_{w1}, \dots, i_{wd})$. We will write $(\ui,\uj) \sim (\uk,\ul)$
if $(\ui,\uj)$ and $(\uk,\ul)$ lie in the same
orbit for the associated diagonal action of $S_d$
on $I(m|n,d) \times I(m|n,d)$.
We say that a double index $(\ui,\uj) \in I(m|n,d) \times I(m|n,d)$
is {\em strict} if
$(\bar{i_r} + \bar{j_r})(\bar{i_s}+\bar{j_s}) = \0$
whenever $(i_r,j_r) = (i_s,j_s)$ for $1 \leq r < s \leq d$.
Let $I^2(m|n,d)$ denote the set of all strict double indexes.
Note $(\ui,\uj)$ is strict if and only if the element
$$
\tilde T_{\ui,\uj} := \tilde T_{i_1,j_1} \cdots \tilde T_{i_d,j_d} \in A(m|n,d)
$$
is non-zero. 
Moreover, if $\Omega(m|n,d)$ is a fixed set of orbit representatives
for the action of $S_d$ on $I^2(m|n,d)$, then
the elements 
$\{\tilde T_{\ui,\uj}\}_{(\ui,\uj) \in \Omega(m|n,d)}$ give a basis 
for $A(m|n,d)$.
Given $(\ui,\uj), (\uk,\ul) \in I^2(m|n,d)$ with
$(\ui,\uj) \sim (\uk,\ul)$, we define a sign
$\sigma(\ui,\uj;\uk,\ul)$ by
\begin{equation}
\sigma(\ui,\uj;\uk,\ul) = \gamma(\epsilon_{\ui} + \epsilon_{\uj}, w)
\end{equation}
if $w$ is any element of $S_d$ with $(\ui,\uj) \cdot w = (\uk,\ul)$.
Note $\tilde T_{\uk,\ul} = \sigma(\ui,\uj;\uk,\ul) \tilde T_{\ui,\uj}$.

For $(\ui,\uj) \in I^2(m|n,d)$, let
$\xi_{\ui,\uj} \in S(m|n,d)$ be the unique element satisfying
$$
\xi_{\ui,\uj}(\tilde T_{\ui,\uj}) = 
\alpha(\epsilon_\ui + \epsilon_{\uj}, \epsilon_{\ui} + \epsilon_{\uj}),\qquad
\xi_{\ui,\uj} (\tilde T_{\uk,\ul}) = 0
\hbox{ for all $(\uk,\ul) \not\sim (\ui,\uj)$}.
$$
The elements
$\{\xi_{\ui,\uj}\}_{(\ui,\uj) \in \Omega(m|n,d)}$ give a basis
for $S(m|n,d)$.
Given $\ui \in I(m|n,d)$, let
$v_{\ui} = v_{i_1} \otimes\dots\otimes v_{i_d} \in V^{\otimes d},$
giving us a basis 
$\{v_{\ui}\}_{\ui \in I(m|n,d)}$ for the tensor space $V^{\otimes d}$.
Since $V^{\otimes d}$ is a polynomial representation of degree $d$,
there is a naturally induced representation 
\begin{equation}
\rho_d:S(m|n,d) \rightarrow \End_k(V^{\otimes d}).
\end{equation}
Also let $e_{\ui,\uj} = e_{i_1,j_1} \otimes \dots \otimes e_{i_d,j_d}
\in \End_k(V)^{\otimes d} \cong \End_k(V^{\otimes d})$, so 
$e_{\ui,\uj} v_{\uk} = \delta_{\uj,\uk}
\alpha(\epsilon_{\ui} + \epsilon_{\uj}, \epsilon_{\uk})
v_{\ui}$ for $\ui,\uj,\uk \in I(m|n,d)$.

\begin{lemma}\label{prodform}
The representation $\rho_d:S(m|n,d) \rightarrow \End_k(V^{\otimes d})$
is faithful and satisfies
$$
\rho_d(\xi_{\ui,\uj}) = \sum_{(\uk,\ul) \sim (\ui,\uj)}
\sigma(\ui,\uj;\uk,\ul) e_{\uk,\ul}
$$
for each $(\ui,\uj) \in I^2(m|n,d)$.
Moreover, for $(\ui,\uj), (\uk,\ul) \in I^2(m|n,d)$,
$$
{\xi}_{\ui, \uj} 
{\xi}_{\uk, \ul} = \sum_{(\us, \ut) \in \Omega(m|n, d)} a_{\ui,\uj,\uk, \ul, \us, \ut} {\xi}_{\us, \ut}
$$
where
$a_{\ui,\uj,\uk, \ul, \us, \ut} = \sum
\sigma(\ui,\uj;\us,\uh)\sigma(\uk,\ul;\uh,\ut)
\alpha(\epsilon_\uh + \epsilon_\ut, \epsilon_\us + \epsilon_\uh)$,
summing over all
$\uh \in I(m|n,d)$ with 
$(\us,\uh) \sim (\ui,\uj), (\uh,\ut) \sim (\uk,\ul)$.
\end{lemma}

\begin{proof}
We first observe that the structure map
$\eta:V^{\otimes d} \rightarrow V^{\otimes d} \otimes A(m|n,d)$ satisfies
$$
\eta(v_{\uj}) = \sum_{\ui \in I(m|n,d)}
(-1)^{\bar\ui (\bar\ui + \bar\uj)}
\alpha(\epsilon_{\ui} + \epsilon_{\uj}, \epsilon_{\ui})
v_{\ui} \otimes \tilde T_{\ui,\uj},
$$
where $\bar\ui = \bar i_1+\dots+\bar i_d, \bar\uj = 
\bar j_1+\dots + \bar j_d$.
Using this we calculate from the definition of the action:
\begin{align*}
\xi_{\ui,\uj} v_{\ul} &=
(1 \bar\otimes \xi_{\ui,\uj}) \sum_{\uk \in I(m|n,d)}\!\!\!
(-1)^{\bar\uk (\bar\uk + \bar\ul)}
\alpha(\epsilon_{\uk} + \epsilon_{\ul}, \epsilon_{\uk})
v_{\uk} \otimes \tilde T_{\uk,\ul}\\
&=
\!\!\!\sum_{\uk \in I(m|n,d)}\!\!\!
\alpha(\epsilon_{\uk} + \epsilon_{\ul}, \epsilon_{\uk})
\xi_{\ui,\uj}(\tilde T_{\uk,\ul})
v_{\uk} \\
&= \!\!\!\sum_{(\uk,\ul) \sim (\ui,\uj)}\!\!\!
\sigma(\ui,\uj;\uk,\ul)
\alpha(\epsilon_{\uk}+\epsilon_{\ul}, \epsilon_{\ul})
v_\uk.
\end{align*}
This prove the formula for $\rho_d(\xi_{\ui,\uj})$.
Hence, $\rho_d$ is injective since
the elements $\{\rho_d(\xi_{\ui,\uj})\}_{(\ui,\uj) \in \Omega(m|n,d)}$ 
are clearly linearly independent.
Finally, to derive the product rule, note that
$e_{\ui,\uj} e_{\uk,\ul}
= \delta_{\uj,\uk} \alpha(\epsilon_{\ui}+\epsilon_{\uj},
\epsilon_{\uk} + \epsilon_{\ul}) e_{\ui,\ul}$.
Using this it is easy to compute the product
$\xi_{\ui,\uj} \xi_{\uk,\ul}$ working in the
representation $\rho_d$.
\end{proof}

We next define a right action of the
symmetric group $S_d$ on 
$V^{\otimes d}$ by letting
\begin{equation}\label{nra}
v_{\ui} (j\:\:j+1)
= 
(-1)^{\bar{i}_{j}\bar{i}_{j+1}}
v_{i_{1}}\otimes \dotsb \otimes v_{i_{j+1}} \otimes v_{i_{j}} \otimes \dotsb \otimes v_{i_{d}}
\end{equation}
for each $\ui = (i_1,\dots,i_d) \in I(m|n,d)$ and each $1 \leq j < d$.
For arbitrary $w \in S_d$, we have that
\begin{equation}\label{geq}
v_{\ui} w = \gamma(\epsilon_{\ui}, w) v_{\ui \cdot w}.
\end{equation}
Note the right action of $S_d$ 
is by even $G$-supermodule automorphisms, so it automatically commutes
with the
left action of $S(m|n,d)$ on $V^{\otimes d}$.
The following theorem is well-known,
see for example \cite{BR1, S2}.

\begin{theorem}
$\rho_d:S(m|n,d) {\rightarrow} \End_{kS_d} (V^{\otimes d})$
is an isomorphism.
\end{theorem}

\begin{proof}
We have already shown in Lemma~\ref{prodform}
that $\rho_d$ is injective and that it maps $S(m|n,d)$ into
$\End_{kS_d}(V^{\otimes d})$. For surjectivity, 
take an arbitrary $\theta:V^{\otimes d} \rightarrow V^{\otimes d}$
commuting with the right action of $S_d$. Write
$$
\theta = \sum_{\ui,\uj \in I(m|n,d)} a_{\ui,\uj} e_{\ui,\uj}
$$
for some coefficients $a_{\ui,\uj}$.
Since $\theta$ commutes with each $w \in S_d$, we have that
$(\theta v_{\uj}) w = \theta(v_{\uj} w)$.
A computation using (\ref{geq}) gives that
$$
\gamma(\epsilon_{\uj}, w)
\alpha(\epsilon_{\ui \cdot w} + \epsilon_{\uj \cdot w}, 
\epsilon_{\uj \cdot w})
a_{\ui \cdot w, \uj \cdot w}
= 
\alpha(\epsilon_{\ui} + \epsilon_{\uj}, \epsilon_{\uj})
\gamma(\epsilon_{\ui}, w)
a_{\ui,\uj}
$$
for each $\ui,\uj$. Simplifying this using the definitions of
$\alpha$ and $\gamma$ then gives that
$$
a_{\ui \cdot w, \uj \cdot w} = a_{\ui, \uj} 
\gamma(\epsilon_{\ui} +
\epsilon_{\uj}, w)
$$
Now note that if $(\ui,\uj)$ is not strict, 
we can choose a transposition
$w \in S_d$ so that $\ui \cdot w = \ui, \uj \cdot w = \uj$
and $\gamma(\epsilon_{\ui} +
\epsilon_{\uj}, w) = -1$. Hence, $a_{\ui,\uj} = 0$ in that case.
Otherwise, if $(\ui,\uj)$ is strict and $(\uk,\ul) \sim (\ui,\uj)$, 
we have shown that
$a_{\uk,\ul} = \sigma(\ui,\uj;\uk,\ul) a_{\ui,\uj}$.
It follows easily that $\theta$ is a linear combination of the
elements $\rho_d(\xi_{\ui,\uj})$ computed in Lemma~\ref{prodform}.
Hence, $\rho_d$ is onto.
\end{proof}

We call a weight $\lambda = \sum_{i=1}^{m+n} \lambda_i \varepsilon_i
\in X(T)$ a {\em polynomial weight} if $\lambda_i \geq 0$ for all $i$.
Let $\Lambda(m|n,d)$ denote the set of all 
such polynomial weights satisfying in addition
$\lambda_1+\dots+\lambda_{m+n} = d$.
Note this is exactly the set of weights arising in the $G$-supermodule
$V^{\otimes d}$. 
For $\ui 
\in I(m|n,d)$, let $\wt(\ui) \in \Lambda(m|n,d)$ 
denote the weight of the vector $v_{\ui}$,
so $\wt(\ui) = \sum_{i=1}^{m+n} \lambda_i \varepsilon_i$ where
there are $\lambda_1$ $1$'s, $\lambda_2$ $2$'s, \dots
appearing in the tuple $(i_1,\dots,i_d)$.
Conversely, given $\lambda \in \Lambda(m|n,d)$, let $\ui_\lambda$
denote the 
tuple $(1,\dots,1,2,\dots,2,\dots) 
\in I(m|n,d)$ with $\lambda_i$ $i$'s for each $i$.
Let
\begin{equation}\label{wtidem}
\xi_{\lambda}= \xi_{\ui_\lambda, \ui_\lambda} \in S(m|n,d).
\end{equation}
We note 
that if $M$ is a polynomial representation of $G$ of 
degree $d$ then, by the argument in \cite[3.2]{G1}, the subspace
$\xi_\lambda M$ is exactly the $\lambda$-weight space $M_\lambda$ of $M$
as defined in (\ref{wtsp}).
An elementary calculation using the 
product rule from Lemma~\ref{prodform} shows:

\begin{lemma}\label{L:weightidempotents}
For $(\ui,\uj) \in I^{2}(m|n,d)$,
\begin{align*} 
\xi_{\lambda}\xi_{\ui,\uj}
&=                           
\begin{cases} \xi_{\ui,\uj} & \mbox{if } \wt(\ui)=\lambda, \\
                                          0                   & \mbox{otherwise, }
                          \end{cases} & 
\xi_{\ui,\uj}\xi_{\lambda} &= 
                          \begin{cases} \xi_{\ui,\uj} & \mbox{if } 
\wt(\uj)=\lambda, \\
0                   & \mbox{otherwise.} 
                          \end{cases}
\end{align*}
In particular, $\{\xi_{\lambda}\}_{\lambda \in \Lambda(m|n,d)}$ is a set of 
mutally orthogonal even idempotents whose sum is the identity in $S(m|n,d)$.
\end{lemma}

Now we turn to the problem of classifying the irreducible
polynomial representations of $G$, or equivalently, the irreducible
$S(m|n,d)$-supermodules for all $d \geq 0$.
Suppose for
some $\lambda \in X^+(T)$ that 
the irreducible $G$-supermodule
$L(\lambda)$ is a polynomial representation of degree $d$. 
Since {\em all} its weights are polynomial weights, $\lambda$ must 
belong to the set
$$
\Lambda^+(m|n,d) :=
\{\lambda \in \Lambda(m|n,d) \:|\:\lambda_1 \geq \dots \geq \lambda_m, \lambda_{m+1}\geq \dots \geq \lambda_{m+n}\}
$$
of dominant polynomial weights of degree $d$. However,
$$
\Lambda^{++}(m|n,d) := \{\lambda \in \Lambda^+(m|n,d)\:|\:L(\lambda)
\hbox{ is a polynomial representation}\}
$$
will in general be a proper subset of $\Lambda^+(m|n,d)$, unlike the
purely even case.
The following lemma is an immediate consequence of Theorem~\ref{hwt}
and the general remarks made at the beginning
of the section.

\begin{lemma}\label{ic} The supermodules
$\{L(\lambda)\}_{\lambda \in \Lambda^{++}(m|n,d)}$
form a complete set of pairwise non-isomorphic irreducible
$S(m|n,d)$-supermodules.
\end{lemma}

Of course the main problem now is to describe the set
$\Lambda^{++}(m|n,d)$ combinatorially! 
Over fields of characteristic $0$, the answer is well-known,
see \cite{BR1} or \cite{S2}.
In positive characteristic, Donkin has given a
combinatorial description of $\Lambda^{++}(m|n,d)$ under the assumption
that $d \leq m$, see  \cite[2.3(4)]{D1}.
We give here an alternative proof of Donkin's result, and
describe $\Lambda^{++}(m|n,d)$ in general in Theorem~\ref{last} in the 
next section.

\begin{theorem}\label{dt}
Assume $d \leq m$. Then,
$$
\Lambda^{++}(m|n,d) = \{\lambda \in \Lambda^+(m|n,d)
\:|\:\lambda_{m+1} \equiv \dots \equiv \lambda_{m+n}
\equiv 0 \!\!\!\pmod{p}\}.
$$
\end{theorem}

\begin{proof}
We first recall the argument of Donkin from \cite[2.3(4)]{D1}
to show that all $\lambda \in \Lambda^+(m|n,d)$
with $\lambda_{m+1} \equiv \dots \equiv \lambda_{m+n} \equiv 0 \pmod{p}$
belong to $\Lambda^{++}(m|n,d)$, i.e. that $L(\lambda)$ is a polynomial
representation for all such $\lambda$.
Let $\omega_i = \varepsilon_1+\dots+\varepsilon_i$.
For any $r_1,\dots,r_m \geq 0$, the polynomial representation
$$
V^{\otimes r_1}
\otimes 
\left({\bigwedge}^2 V\right)^{\otimes r_2}
\otimes\dots\otimes
\left({\bigwedge}^m V\right)^{\otimes r_m}
$$
has unique highest weight $r_1 \omega_1 + \dots + r_m \omega_m$.
Hence $L(r_1 \omega_1 + \dots + r_m \omega_m)$ is a composition factor
of a polynomial representation, so polynomial.
Now given an arbitrary 
$\lambda \in \Lambda^{++}(m|n,d)$
with $\lambda_{m+1} \equiv \dots \equiv \lambda_{m+n} \equiv 0 \pmod{p}$,
we can write $\lambda = \mu + p \nu$ where
$\mu \in X^+(T)_{\operatorname{res}}$ is a restricted polynomial weight 
in the sense of Remark~\ref{rks}(ii) satisfying 
$\mu_{m+1} = \dots = \mu_{m+n} = 0$, and $\nu \in X^+(T)$ is an
arbitrary polynomial weight.
Then, the $G$-supermodule
$$
L(\mu) \otimes \Fr^* L_{\ev}(\nu)
$$
has unique
highest weight $\lambda$ so contains $L(\lambda)$ as a composition factor
(actually it equals $L(\lambda)$ by Remark~\ref{rks}(iii), though we do not
need this stronger result).
Since $\mu$ can be expressed in the form 
$r_1 \omega_1 + \dots + r_m \omega_m$ for $r_i \geq 0$, 
$L(\mu)$ is polynomial and 
$L_{\ev}(\nu)$, hence $\Fr^* L_{\ev}(\nu)$, is polynomial
by the classical theory. So $L(\lambda)$ is a polynomial representation, 
and $\lambda \in \Lambda^{++}(m|n,d)$.

Conversely, suppose that
$\lambda \in \Lambda^+(m|n,d)$ has
$\lambda_{m+i} \not\equiv 0 \pmod{p}$ for some $1 \leq i \leq n$.
Pick the minimal such $i$, and observe by the assumption $d \leq m$
that $\lambda_m = 0$.
Let $w \in D_{m,n}$ be the 
permutation $(m\:\:m+1\:\:\dots\:\:m+i)$, so
$$
\Phi_w^+ = \Phi_1^+ - \{\varepsilon_m - \varepsilon_{m+1},\dots,
\varepsilon_m - \varepsilon_{m+i}\}
\cup
\{\varepsilon_{m+1} - \varepsilon_{m},\dots,
\varepsilon_{m+i} - \varepsilon_{m}\}.
$$
Applying Lemma~\ref{st} using the sequence $\varepsilon_m - \varepsilon_{m+1},
\dots, \varepsilon_m - \varepsilon_{m+i}$ of roots, we see that
$L(\lambda) \cong L_w(\lambda - 
\varepsilon_m + \varepsilon_{m+i})$. Thus, $\lambda - \varepsilon_m + \varepsilon_{m+i}$ is a weight of
$L(\lambda)$. Since this is not a polynomial weight, $L(\lambda)$
cannot be a polynomial representation, i.e.
$\lambda \notin \Lambda^{++}(m|n,d)$.
\end{proof}

We now explain how to descend from the Schur superalgebra $S(m|n,d)$
to the symmetric group $S_d$, assuming still that $d \leq m$.
We need the following basic fact about functors
defined by idempotents, cf. \cite[6.2]{G1} 
or \cite[Corollary 2.13]{BK1}. Recall that if $\xi \in S(m|n,d)$
is an even idempotent and $M$ is an $S(m|n,d)$-supermodule, we can
view $\xi M$ naturally as a supermodule over the
subalgebra $\xi S(m|n,d) \xi$ of $S(m|n,d)$.

\begin{lemma}\label{il}
Let $\xi \in S(m|n,d)$ be an even idempotent.
For $\lambda \in \Lambda^{++}(m|n,d)$, $\xi L(\lambda)$ is either zero
or it is an irreducible $\xi S(m|n,d) \xi$-supermodule.
Moreover, the non-zero $\xi L(\lambda)$'s give a complete set of
pairwise non-isomorphic irreducible $\xi S(m|n,d) \xi$-supermodules.
\end{lemma}

One checks the following lemma using the product rule
from Lemma~\ref{prodform}.

\begin{lemma}\label{littlecheck}
Assume $d \leq m$ and let $\omega = \sum_{i=1}^d \varepsilon_i
\in \Lambda(m|n,d)$.
Then, for any $\ui \in I(m|n,d)$ and any $x \in S_d$,
$
\xi_{\ui,\ui_\omega} \xi_{\ui_\omega \cdot x, \ui_\omega}
= 
\gamma(\epsilon_{\ui}, x) \xi_{\ui \cdot x, \ui_\omega}.
$
\end{lemma}

Continue with $d\leq m$ and $\omega$ as in Lemma~\ref{littlecheck}.
By Lemma~\ref{L:weightidempotents}, the subalgebra
$\xi_\omega S(m|n,d) \xi_\omega$ of $S(m|n,d)$ has
basis $\{\xi_{\ui_\omega \cdot x, \ui_\omega}\}_{x \in S_d}$.
Lemma~\ref{littlecheck} implies that the map
\begin{equation}\label{e2}
kS_d \rightarrow \xi_\omega S(m|n,d) \xi_\omega,
\qquad 
x \mapsto \xi_{\ui_\omega \cdot x, \ui_\omega}
\end{equation}
is a superalgebra isomorphism. 
>From now on, we will {\em identify} $k S_d$
with the subalgebra 
$\xi_\omega S(m|n,d) \xi_\omega$ of $S(m|n,d)$ in this way.
Then, we can define the {\em Schur functor}
\begin{equation}\label{sf}
f_\omega: S(m|n,d)\hbox{-mod} 
\rightarrow kS_d\hbox{-mod}.
\end{equation}
On an $S(m|n,d)$-supermodule $M$, $f_\omega M$
is the $\omega$-weight space  $\xi_{\omega} M$ of $M$ viewed as a
$kS_d$-supermodule via the identification (\ref{e2}). On a morphism,
the functor $f_\omega$ is defined by restriction.

\begin{remark}\label{morerks}\rm
(i) Recalling the definition of the
action of $S_d$ on $V^{\otimes d}$ from (\ref{geq}), 
Lemma~\ref{littlecheck} also shows that
the map
$V^{\otimes d} \rightarrow S(m|n,d)\xi_{\omega}, \:
v_{\ui} \mapsto \xi_{\ui,\ui_{\omega}}
$
is an isomorphism of $S(m|n,d), kS_d$-bimodules.

(ii) An immediate consequence of (i) is that the Schur functor 
$f_\omega$ can be defined alternatively by
$f_\omega M = \Hom_G(V^{\otimes d}, M)$
for a polynomial $G$-supermodule $M$ of degree $d$, where the $S_d$
action on $f_\omega M$ is induced by the natural right action of
$S_d$ on $V^{\otimes d}$ from (\ref{nra}).

(iii) Another well-known consequence of (i) is the
{\em double centralizer property}: for $d \leq m$,
$\operatorname{End}_{S(m|n,d)} (V^{\otimes d}) = k S_d$.
Indeed, by (i) and properties of idempotents,
$$
\operatorname{End}_{S(m|n,d)} (V^{\otimes d})
\cong
\operatorname{End}_{S(m|n,d)} (S(m|n,d) \xi_{\omega})
\cong \xi_\omega S(m|n,d) \xi_\omega \cong k S_d.
$$
One can extend this result to the case $d \leq m+n$ by similar arguments.
\end{remark}

\vspace{1mm}

Recall now that a {\em partition of $d$} is a sequence
$\lambda = (\lambda_1 \geq \lambda_2 \geq \dots)$ of non-negative integers
satisfying $|\lambda| := \lambda_1+\lambda_2+ \dots = d$.
We usually identify $\lambda$ with its
{\em Young diagram}
$$
\lambda = \{(i,j) \in \mathbb Z_{> 0} \times \mathbb Z_{> 0}\:|\:j \leq \lambda_i\}
$$
and refer to $(i,j) \in \lambda$ as the {\em node} in the $i$th
row and $j$th column.
We say that a partition
$\lambda = 
(\lambda_1 \geq \lambda_2 \geq \dots)$ is 
{\em restricted} 
if either $p = 0$ or $p > 0$ and $\lambda_i - \lambda_{i+1} < p$ for all $i = 1,2,\dots$ (cf. Remark~\ref{rks}(ii)).
Let $\mathscr P(d)$ denote the set of all partitions of $d$, and
$\mathscr{RP}(d) \subseteq \mathscr{P}(d)$ 
denote the set of all {restricted partitions of $d$}.

Assuming still that $d \leq m$, we define an embedding
\begin{equation}\label{x}
x:\mathscr{RP}(d) \hookrightarrow \Lambda^{+}(m|n,d),
\qquad
\lambda \mapsto \sum_{i=1}^m \lambda_i \varepsilon_i.
\end{equation}
By Theorem~\ref{dt}, we actually have that
$x(\lambda) \in \Lambda^{++}(m|n,d)$ for 
$\lambda \in \mathscr{RP}(d)$. So it makes sense to define
\begin{equation}\label{ddef}
D_\lambda = f_\omega L(x(\lambda)).
\end{equation}
The following theorem shows in particular that the
$D_\lambda$'s are non-zero $kS_d$-modules.

\begin{theorem}\label{parm} Assume $d \leq m$.
For $\lambda \in \Lambda^{++}(m|n,d)$,
$f_\omega L(\lambda) \neq 0$ if and only if $\lambda = x(\mu)$
for some $\mu \in \mathscr{RP}(d)$.
Hence, the $kS_d$-modules 
$\{D_\lambda\}_{\lambda \in \mathscr{RP}(d)}$
form a complete set of pairwise non-isomorphic irreducible
$kS_d$-modules.
\end{theorem}

\begin{proof}
Take $\lambda \in \Lambda^{++}(m|n,d)$.
By Theorem~\ref{dt}, if $\lambda \notin x(\mathscr{RP}(d))$, we can
decompose $\lambda = \mu + p \nu$ for polynomial weights
$\mu,\nu \in X^+(T)$
with $\mu$ restricted in the sense of Remark~\ref{rks}(ii)
and with $\nu \neq 0$. But then
$L(\mu) \otimes \Fr^* L(\nu)$
has unique highest weight $\lambda$, so has $L(\lambda)$ as
a composition factor (actually it equals $L(\lambda)$
by Remark~\ref{rks}(iii)).
Since the $\omega$-weight space of
$L(\mu) \otimes \Fr^* L(\nu)$ is zero, this shows that
$f_\omega L(\lambda) = 0$.
Finally, by Lemma~\ref{il}, 
the non-zero $f_\omega L(\lambda)$ with $\lambda \in \Lambda^{++}(m|n,d)$
must give a complete set of pairwise non-isomorphic irreducible $kS_d$-modules.
It is well-known that the number of isomorphism classes of 
the latter is $|\mathscr{RP}(d)|$, hence we must have that
$f_\omega L(\lambda) \neq 0$ for all $\lambda \in x(\mathscr{RP}(d))$.
This completes the proof.
\end{proof}

\begin{remark}\label{lrk}\rm
The $p$-regular partitions from the introduction are the conjugates of
the restricted partitions. However,  
we will work from now on with the parametrization of the irreducible
$kS_d$-modules by restricted partitions according to Theorem~\ref{parm},
though this is not the usual convention made in the literature.
The relationship between our labeling and the standard labeling
of James \cite{J1} is given by
\begin{equation}\label{labrel}
D_\lambda \cong D^{\lambda'} \otimes \sgn.
\end{equation}
One can see this as follows. Embedding $\Lambda(m,d) := \Lambda(m|0,d)$ into
$\Lambda(m|n,d)$ as the set of all weights with $\lambda_{m+1}=\dots=\lambda_{m+n} = 0$, let
$\xi = \sum_{\lambda \in \Lambda(m,d)} \xi_\lambda.$
Then, $\xi S(m|n,d) \xi$ can be identified with the classical
Schur superalgebra $S(m,d) := S(m|0,d)$ of \cite{G1},
see the proof of Theorem~\ref{last} below for a similar construction.
Moreover, given $\lambda \in \Lambda^+(m,d) := \Lambda^+(m|0,d)$, 
$\xi L(\lambda)$ is the irreducible $S(m,d)$-module with highest weight
$\lambda$. Notice that for $\omega$ as in Lemma~\ref{littlecheck},
$\xi_{\omega} \xi = \xi_\omega$. Hence, our Schur functor
$f_\omega$ from representations of $S(m|n,d)$ to representations of
$kS_d$ factors through the
Schur functor in \cite[6.4]{G1} 
from representations of $S(m,d)$ to representations
of $kS_d$.
So \cite[6.4]{G1} implies that $D_\lambda =
 L(x(\lambda)) \cong
D^{\lambda'} \otimes \sgn$ for each $\lambda \in \mathscr{RP}(d)$.
\end{remark}

\section{The Mullineux conjecture}\label{S:Mullineux}

Let $\lambda \in \mathscr{P}(d)$ be a partition of $d$.
The {\em rim} of $\lambda$ is defined to be the set of 
all nodes $(i,j)\in\lambda$ such that $(i+1,j+1)\not\in\lambda$.
The {\em $p$-rim} is a certain subset of 
the rim, defined as the union of the {\em $p$-segments}. The first $p$-segment
is simply the first $p$ nodes of the rim, reading along the rim 
from left to right.
The next $p$-segment is then obtained
by reading off the next $p$ nodes of the rim, but starting 
from the column immediately to the {\em right} of
the rightmost node of the first $p$-segment.
The remaining $p$-segments are obtained by repeating this process. Of course,
all but the last $p$-segment contain exactly $p$ nodes, while the last may
contain less. For example,
let $\lambda=(5,4,3^2,1^2),\ p=5$. The nodes of the $p$-rim 
(which consists of two $p$-segments) are colored in black in the following 
picture.
$$
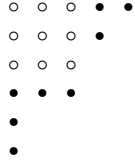
\setlength{\unitlength}{0.007500in}%
\begin{picture}(108,124)(139,514)
\multiput(143,618)(20,0){3}{\circle{5}}
\multiput(143,598)(20,0){3}{\circle{5}}
\multiput(143,578)(20,0){3}{\circle{5}}
\multiput(203,618)(20,0){2}{\circle*{6}}
\multiput(203,598)(20,0){1}{\circle*{6}}
\multiput(143,558)(20,0){3}{\circle*{6}}
\multiput(143,538)(20,0){1}{\circle*{6}}
\multiput(143,518)(20,0){1}{\circle*{6}}
\end{picture}
$$
Let $a(\lambda)$ denote the number of nodes in the $p$-rim
of $\lambda$.

We now define {\em Mullineux conjugation}
$$
\mx:\mathscr{RP}(d) \rightarrow \mathscr{RP}(d),
$$
which is actually
the transpose of the original definition from \cite{M1}
since we are working with restricted rather than regular partitions.
Given $\lambda \in \mathscr{RP}(d)$, 
set $\lambda^{(1)}=\lambda$, and define $\lambda^{(i)}$ to be 
$\lambda^{(i-1)} - \{\text{the $p$-rim of $\lambda^{(i-1)}$}\}$. 
Let $m$ be the largest number such that $\lambda^{(m)}\neq 0$. 
The {\em Mullineux symbol} of $\lambda$ is defined to be the array
$$
G(\lambda) =
\left(
\begin{matrix}
a_1 & a_2 & \dots & a_m \\
r_1 & r_2 & \dots & r_m
\end{matrix}
\right)
$$
where $a_i = a(\lambda^{(i)})$ 
is the number of the nodes in the $p$-rim of $\lambda^{(i)}$ and 
$r_i = \lambda^{(i)}_1$ is the first part of $\lambda^{(i)}$.
The partition $\lambda$
can be uniquely reconstructed from its Mullineux symbol, 
see \cite{M1}.
Now, $\mx(\lambda)$ is defined to be the unique restricted partition with
\begin{equation*}
G(\mx(\lambda)) =
\left(
\begin{matrix}
a_1 & a_2 & \dots & a_m \\
s_1 & s_2 & \dots & s_m
\end{matrix}
\right)
\end{equation*}
where 
\begin{equation}\label{pp}
s_i=\left\{\begin{array}{ll}
a_i-r_i&\hbox{if $a_i \equiv 0 \pmod{p}$,}\\
a_i+1-r_i&\hbox{if $a_i \not\equiv 0 \pmod{p}$.}
\end{array}\right.
\end{equation}
Note in particular that 
the first part of $\mx(\lambda)$ equals $s_1$.

As explained in the introduction, we will be concerned here 
with an equivalent formulation of the Mullineux algorithm
discovered by Xu \cite{Xu1}.
For $\lambda \in \mathscr{P}(d)$, let $\jx(\lambda)$ be the partition
obtained from $\lambda$ 
by deleting every node in the $p$-rim that is at the rightmost end
of a row of $\lambda$ but that is not the $p$th node of a $p$-segment.
Let $j(\lambda) = |\lambda| - |\jx(\lambda)|$ be the total number of nodes
deleted.
For example, with $\lambda = (5,4,3^2,1^2), p = 5$ as above, 
$\jx(\lambda)$ is obtained by deleting the double-circled 
nodes:
$$
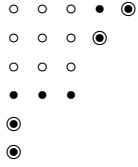
\setlength{\unitlength}{0.007500in}%
\begin{picture}(108,124)(139,514)
\multiput(143,618)(20,0){3}{\circle{5}}
\multiput(143,598)(20,0){3}{\circle{5}}
\multiput(143,578)(20,0){3}{\circle{5}}
\multiput(203,618)(20,0){2}{\circle*{6}}
\multiput(203,598)(20,0){1}{\circle*{6}}
\multiput(143,558)(20,0){3}{\circle*{6}}
\multiput(143,538)(20,0){1}{\circle*{6}}
\multiput(143,518)(20,0){1}{\circle*{6}}
\multiput(143,518)(20,0){1}{\circle{10}}
\multiput(203,598)(20,0){1}{\circle{10}}
\multiput(223,618)(20,0){1}{\circle{10}}
\multiput(143,538)(20,0){1}{\circle{10}}
\end{picture}
$$
Hence $j(\lambda) = 4$.
Note the definitions of the maps $\jx$ and $j$
make sense for arbitrary partitions,
though to prove the Mullineux conjecture we only need to apply them to
restricted partitions.
In general, one has that $\jx(\mu+p\nu) = \jx(\mu)+p\nu$ and
$j(\mu+p\nu)= j(\mu)$. 

Recalling that $a(\lambda)$ is the number of nodes in the $p$-rim of $\lambda$,
we note for arbitrary $\lambda \in \mathscr{P}(d)$ that
\begin{equation}\label{sd}
j(\lambda) = \left\{\begin{array}{ll}
a(\lambda)-\lambda_1&\hbox{if $a(\lambda) \equiv 0 \pmod{p}$,}\\
a(\lambda)+1-\lambda_1&\hbox{if $a(\lambda) \not\equiv 0 \pmod{p}$.}\\
\end{array}\right.
\end{equation}
Comparing with (\ref{pp}), this shows that for restricted $\lambda$, 
$j(\lambda)$ is 
the first part of the partition $\mx(\lambda)$.
More generally, it is proved in \cite[Proposition 3.4]{BOX} that
for restricted $\lambda$,
$\mx(\jx(\lambda)) = \rx(\mx(\lambda))$, 
where $\rx$ denotes first row removal.
Using this fundamental 
fact, the following theorem of Xu \cite{Xu1} 
follows easily:

\begin{theorem}\label{xua}
For $\lambda \in \mathscr{RP}(d)$, 
$\mx(\lambda)$ is the partition $\mu$ with $\mu_i = j(\jx^{i-1}(\lambda))$.
\end{theorem}

We will refer to the algorithm for computing $\mx(\lambda)$ given by
Theorem~\ref{xua} as Xu's algorithm. 
For an example, take $\lambda = (5,4, 3^2, 1^2), p = 5$ as above.
Then $\jx(\lambda) = (4,3^3),
\jx^2(\lambda) = (3^2, 2^2), 
\jx^3(\lambda) = (3^2,1^2), \jx^4(\lambda) = (3^2),
\jx^5(\lambda) = (2^2), \jx^6(\lambda) = (1^2), \jx^7(\lambda) = 0$.
Hence, $\mx(\lambda) = (4,3,2^5)$.

We next explain the relationship
between Xu's algorithm and Serganova's algorithm from 
Theorem~\ref{salg}. 
The main step is to prove the following alternative description
of the map $\jx$.

\begin{lemma}\label{jj}
Suppose that $\lambda  \in \mathscr{P}(d)$ with
$\lambda_{m+1} = 0$.
Define $x_1,x_2,\dots \in \{0,1\}$ 
by $x_{m+1} = x_{m+2} = \dots = 0$ and 
$$
x_i = \left\{ 
\begin{array}{ll}
1 &\hbox{if $\lambda_i + x_{i+1} + x_{i+2} + \dots \not\equiv 0 \pmod{p}$,}\\
0 &\hbox{if $\lambda_i + x_{i+1} + x_{i+2} + \dots \equiv 0 \pmod{p}$,}
\end{array}\right.
$$
for $i=m,m-1,\dots,1$.
Then, $\jx(\lambda)$ is the partition $\mu$ with $\mu_i = \lambda_i - x_i$.
\end{lemma}

\begin{proof}
Proceed by induction on $m$, the case $m = 0$ being vacuous.
For the induction step, take $\lambda \in \mathscr{RP}(d)$ 
with $\lambda_{m+1} = 0$.
Define $x_1,x_2,\dots$ and $\mu$ according to the statement of the lemma.
By the induction hypothesis, we get that
$\jx(\rx(\lambda)) = \rx(\mu)$, which shows in particular that
$j(\rx(\lambda)) = |\rx(\lambda)| - |\rx(\mu)| = 
x_2 + x_3 + \dots$.
To complete the proof, it remains to show that 
the first part of $\jx(\lambda)$
is equal to $\lambda_1 - x_1$, or equivalently,
$j(\lambda) = x_1 + x_2 + x_3 + \dots$.

If $a(\rx(\lambda))
\equiv 0 \pmod{p}$, then
all the $p$-segments in the $p$-rim of $\rx(\lambda)$ have $p$ nodes in
them. Hence, the node $(1, \lambda_2)$ does not belong to the
$p$-rim of $\lambda$. Using (\ref{sd}) for the second equality, 
we therefore get that
$$
a(\lambda) = 
\lambda_1 - \lambda_2 + 
a(\rx(\lambda))
=
\lambda_1 + j(\rx(\lambda)).
$$
Otherwise, if $a(\rx(\lambda)) \not\equiv 0 \pmod{p}$, then
the last $p$-segment of $\rx(\lambda)$ has less than $p$ nodes in it.
This implies that
the node $(1, \lambda_2)$ must belong to the $p$-rim of $\lambda$,
so
$$
a(\lambda) = \lambda_1 - \lambda_2 + 1 + a(\rx(\lambda))
=
\lambda_1 + j(\rx(\lambda)).
$$
Thus in either case, we have shown that
$$
a(\lambda) = \lambda_1 + j (\rx(\lambda)) = 
 \lambda_1 + x_2 + x_3 + \dots.
$$
If this is zero mod $p$,
then $x_1 = 0$ and $j(\lambda) = a(\lambda) - \lambda_1$ by
(\ref{sd}).
If it is non-zero mod $p$, then
$x_1 = 1$ and $j(\lambda) = a(\lambda) + 1 - \lambda_1$.
Either way, $j(\lambda) = x_1 + x_2 + \dots$ as required.
\end{proof}

Assume now that $m,n \geq d$.
Recall the definition of the embedding
$x:\mathscr{RP}(d) \hookrightarrow \Lambda^+(m|n,d)$
from (\ref{x}). Instead, define
\begin{equation}\label{y}
y:\mathscr{RP}(d) \hookrightarrow \Lambda^+(m|n,d),
\qquad
\lambda \mapsto \sum_{i=1}^{n} \lambda_i \varepsilon_{m+i}.
\end{equation}
Let $\sim:X^+(T) \rightarrow X^+(T)$ be the bijection defined
combinatorially according to Theorem~\ref{salg}.
Then:

\begin{lemma}\label{cpart}
For $m,n \geq d$ and $\lambda \in \mathscr{RP}(d)$, 
$\widetilde{x(\lambda)} = y(\mx(\lambda))$.
\end{lemma}

\begin{proof}
Compute $\widetilde{x(\lambda)}$ using Theorem~\ref{salg} and the ordering
$$
\varepsilon_m - \varepsilon_{m+1}, \dots, \varepsilon_1 - \varepsilon_{m+1};
\varepsilon_m - \varepsilon_{m+2}, \dots, \varepsilon_1 - \varepsilon_{m+2};
\dots;
\varepsilon_m - \varepsilon_{m+n}, \dots, \varepsilon_1 - \varepsilon_{m+n}.
$$
After the first $m$ steps of the process, $x(\lambda)$
has been replaced by 
$x(\lambda) - \sum_{i=1}^m x_i \varepsilon_i 
+ j(\lambda) \varepsilon_{m+1}$, where $x_1,\dots,x_m$ are defined
as in Lemma~\ref{jj}.
Lemma~\ref{jj} shows this is exactly the weight
$x(\jx(\lambda)) + j(\lambda) \varepsilon_{m+1}.$
Repeating the argument for the next $m$ steps of Serganova's
algorithm, we see similarly that
$x(\jx(\lambda)) + j(\lambda) \varepsilon_{m+1}$
gets replaced by the weight
$x(\jx(\jx(\lambda))) + j(\lambda) \varepsilon_{m+1} + j
(\jx(\lambda)) \varepsilon_{m+2}$.
Continuing in this way and using Xu's
Theorem~\ref{xua}, we get after the $nm$th step of Serganova's 
algorithm that
$\widetilde{x(\lambda)} = y(\mx(\lambda))$.
\end{proof}

At last we are ready to prove the Mullineux conjecture,
see also (\ref{labrel}).

\begin{theorem}
For $\lambda \in \mathscr{RP}(d)$,
$D_\lambda \otimes \sgn \cong D_{\mx(\lambda)}$.
\end{theorem}

\begin{proof}
Take $m = n \geq d$.
Let $\sigma:G \rightarrow G$ be the supergroup automorphism 
defined for a commutative superalgebra $A$ and a matrix
$g \in G(A)$ of the form (\ref{E:matrix}) by
\begin{equation*} 
\left(
\begin{array}{l|l}
W&X\\\hline
Y&Z
\end{array}
 \right) \mapsto
\left(\begin{array}{l|l}
Z&Y\\\hline
X&W
\end{array}
 \right).
\end{equation*}
Given any $G$-supermodule $M$, we let $\sigma^* M$ denote the
$G$-supermodule equal to $M$ as a vector superspace, but with new
action defined by $g \cdot m = \sigma(g) m$
for all $g \in G(A), m \in M \otimes A$ and all commutative superalgebras
$A$.
In particular, $\sigma^*(V^{\otimes d})$ denotes the tensor space
$V^{\otimes d}$ with the action of $G$ twisted by $\sigma$ and
with the original $S_d$-action from (\ref{nra}).
We also write $V^{\otimes d} \otimes \sgn$ for the $G$-supermodule
$V^{\otimes d}$ but with the action of $S_d$ twisted by tensoring with
$\sgn$.

Let $\sigma:\{1,\dots,2n\} \rightarrow \{1,\dots,2n\}$ be the
map $i \mapsto i+n$ if $i \leq n$, $i \mapsto i - n$ if $i \geq n+1$.
Extend $\sigma$ to a map $\sigma:I(n|n,d) \rightarrow I(n|n,d)$
mapping $\ui = (i_1,\dots,i_d)$ to $\sigma(\ui) =
(\sigma(i_1), \dots, \sigma(i_d))$.
Define a map $\sigma:V \rightarrow \sigma^* V, v_i \mapsto v_{\sigma(i)}$.
Obviously, this is an odd isomorphism of $G$-supermodules.
Hence, the map
$$
\sigma^{\otimes d}:V^{\otimes d} \otimes \sgn  
\rightarrow \sigma^* (V^{\otimes d}),\quad
v_{\ui} \otimes 1 \mapsto
 (-1)^{(d-1) \bar i_1 + (d-2) \bar i_2 + 
\dots + \bar i_{d-1}}
v_{\sigma(\ui)}
$$
is an isomorphism of $G$-supermodules.
Using this formula, it is easy to check that 
the map $\sigma^{\otimes d}$ commutes with the action of a 
simple transposition $(j\:\:j+1) \in S_d$.
Hence, $\sigma^{\otimes d}$ is an isomorphism of $G, S_d$-bimodules.
It follows immediately that for any $G$-supermodule $M$,
there is a natural isomorphism
$$
\Hom_G(V^{\otimes d}, \sigma^* M) =
\Hom_G(\sigma^*(V^{\otimes d}), M)
\cong
\Hom_G(V^{\otimes d} \otimes \sgn, M)
$$
of $kS_d$-modules. 
Hence, recalling Remark~\ref{morerks}(ii), we have 
a natural isomorphism
\begin{equation}\label{mainpoint}
(f_\omega M) \otimes \sgn\cong
f_{\omega} (\sigma^* M)
\end{equation}
of $kS_d$-modules for any $S(n|n,d)$-supermodule $M$.

We now apply (\ref{mainpoint}) to $M = L(x(\lambda))$.
By (\ref{Tdef})  and Lemma~\ref{cpart},
$$
L(x(\lambda)) \cong L_{w_1}(\widetilde{x(\lambda)})
\cong L_{w_1}(y(\mx(\lambda))).
$$
The automorphism $\sigma$ of $G$ swaps the Borel subgroups $B_1$
and $B_{w_1}$ and interchanges the two diagonal blocks in the torus $T$. 
Hence,
$$
\sigma^* L_{w_1}(y(\mx(\lambda))) \cong L(x(\mx(\lambda))).
$$
So by (\ref{ddef}) and (\ref{mainpoint}), we get
$$
D_\lambda \otimes \sgn
=
(f_\omega L(x(\lambda))) \otimes \sgn \cong
f_\omega (\sigma^* L(x(\lambda)))
\cong f_\omega L(x(\mx(\lambda))) = D_{\mx(\lambda)}.
$$
This completes the proof.
\end{proof}

We conclude the article by completing the combinatorial description
of the set $\Lambda^{++}(m|n,d)$ that parametrizes the irreducible polynomial
representations of $G$ of degree $d$ in Lemma~\ref{ic}. For $\lambda 
\in \Lambda^+(m|n,d)$, we will use the notation
$t(\lambda)$ for the partition
$(\lambda_{m+1}, \lambda_{m+2}, \dots, \lambda_{m+n})$, i.e. 
the ``tail'' of $\lambda$.
Also recall the definition of $j$ from (\ref{sd}).

\begin{theorem}\label{last} For arbitrary $m,n,d$, we have that
\begin{align*}
\Lambda^{++}(m|n,d) &= \{\lambda \in \Lambda^+(m|n,d)\:|\:
j(t(\lambda)) \leq \lambda_m\}.
\end{align*}
\end{theorem}

\begin{proof}
Pick $M \geq m, N \geq n$ such that $M \geq d$ and $N = M-m$. 
Throughout the proof, 
we will identify $GL(m|n)$ with the closed subgroup of $GL(M|N)$
consisting (for each commutative superalgebra $A$) 
of all invertible matrices of
the form
\begin{equation*}
\left(
\begin{array}{ll|ll}
W&0&X&0\\
0&I_{M-m}&0&0\\
\hline
Y&0&Z&0\\
0&0&0&I_{N-n}
\end{array}
 \right),
\end{equation*}
where $W, X, Y, Z$ are as in (\ref{E:matrix}).
Embed $\Lambda(m|n,d)$ (resp. $\Lambda^+(m|n,d)$)
into $\Lambda(M|N,d)$ (resp. $\Lambda^+(M|N,d)$)
as the set of all $\lambda$ 
with $\lambda_{m+1} = \dots = \lambda_M
= \lambda_{M+n+1} = \dots = \lambda_{M+N} = 0$, and 
embed $I(m|n,d)$ into $I(M|N,d)$ as the 
set of all $d$-tuples $\ui$ with 
entries belonging to the set $\{1,\dots,m, M+1,\dots,M+n\}$.
Let
$$
\xi = \sum_{\lambda \in \Lambda(m|n,d)}
\xi_\lambda \in S(M|N,d).
$$
The embedding $GL(m|n) \hookrightarrow GL(M|N)$ induces an isomorphism
between the Schur superalgebra $S(m|n,d)$ and the subalgebra
$\xi S(M|N,d) \xi$ of $S(M|N,d)$.
The element $\xi_{\ui,\uj}$ of $S(m|n,d)$ for
$\ui,\uj \in I(m|n,d)$ corresponds to the element of $\xi S(M|N,d)\xi$ 
with the same name. 
We will denote the irreducible 
$S(M|N,d)$-supermodule of highest weight
$\lambda \in \Lambda^{++}(M|N,d)$ by $L(\lambda)$, and 
the irreducible 
$S(m|n,d)$-supermodule of highest weight
$\lambda \in \Lambda^{++}(m|n,d)$ by $L'(\lambda)$.

Let $w \in D_{M,N}$ be the permutation
\begin{equation*}
w=\left( 
\begin{matrix} m+1 & m+2 &\dots& M \\
               M+1 & M+2 & \dots& M+N
\end{matrix} \right)
\end{equation*}
This defines a Borel subgroup $B_w$ of $G = GL(M|N)$, 
a set $\Phi_w^+$ of positive roots
and a dominance ordering $\leq_w$ on $X(T)$. Explicitly, for
a commutative superalgebra $A$, $B_w(A)$ consists of all matrices
in $G(A)$ of the form
$$
\left(
\begin{array}{ll|l}
P&X&Y\\
0&Q&0\\
\hline
0&Z&R
\end{array}
 \right),
$$
where $P$ is an upper triangular $m \times m$ matrix,
$Q, R$ are upper triangular $N \times N$ matrices, and $X, Y, Z$ are arbitrary.
Like in (\ref{Tdef}), 
define a bijection $r:X(T) \rightarrow X(T)$ by the rule
$$
L_1(\lambda) \cong L_w(r(\lambda)).
$$
The key observation is that if $\lambda$ is a weight $\notin \Lambda(m|n,d)$,
then every $\mu \leq_w \lambda$ is also
$\notin \Lambda(m|n,d)$.
Hence, since the idempotent $\xi$ is just projection onto the weight
spaces belonging to $\Lambda(m|n,d)$, we see that for
$\lambda \in \Lambda^{++}(M|N,d)$,
$\xi L(\lambda) \cong \xi L_w(r(\lambda))$ is non-zero
if and only if $r(\lambda) \in \Lambda^+(m|n,d)$.
Moreover, in that case, $r(\lambda)$ is the highest weight of
$\xi L(\lambda)$ with respect to the standard dominance ordering on
$\Lambda(m|n,d)$. Viewing $\xi L(\lambda)$ as an $S(m|n,d)$-supermodule
via the identification $S(m|n,d) = \xi S(M|N,d) \xi$, 
we have shown:
\begin{equation}
\xi L(\lambda) \cong \left\{
\begin{array}{ll}
L'(r(\lambda))&\hbox{if $r(\lambda) \in \Lambda^+(m|n,d)$,}\\
0&\hbox{otherwise.}
\end{array}\right.
\end{equation}
Invoking Lemma~\ref{il}, this means that 
\begin{align*}
\Lambda^{++}(m|n,d) 
&= 
r(\Lambda^{++}(M|N,d))
\cap \Lambda^+(m|n,d)\\
&= \{\lambda \in \Lambda^+(m|n,d)\:|\:
r^{-1}(\lambda) \in \Lambda^{++}(M|N,d)\}.
\end{align*}

Now we compute $r^{-1}(\lambda)$ for $\lambda \in \Lambda^{+}(m|n,d)
\subseteq \Lambda^+(M|N,d)$.
Let $t(\lambda)$ denote the partition $(\lambda_{M+1}, \lambda_{M+2},\dots,
\lambda_{M+N})$, i.e. the tail of $\lambda$ as in the statement of the theorem.
Write $t(\lambda) = \mu + p \nu$ for partitions $\mu,\nu$ with 
$\mu$ restricted, so
$$
\lambda = 
\sum_{i=1}^m \lambda_i \varepsilon_i
+ \sum_{i=1}^N \mu_i \varepsilon_{M+i}
+ p \sum_{i=1}^N \nu_i \varepsilon_{M+i}.
$$
Applying Lemma~\ref{st} repeatedly to the root sequence
\begin{multline*}
\varepsilon_{M+N} - \varepsilon_{m+1},
\dots,
\varepsilon_{M+1} - \varepsilon_{m+1};
\varepsilon_{M+N} - \varepsilon_{m+2},
\dots,
\varepsilon_{M+1} - \varepsilon_{m+2};
\dots;\\
\varepsilon_{M+N} - \varepsilon_{M},
\dots,
\varepsilon_{M+1} - \varepsilon_{M}.
\end{multline*}
and arguing as in the proof of
Lemma~\ref{cpart}, one gets that
$$
r^{-1}(\lambda) = 
\sum_{i=1}^m \lambda_i \varepsilon_i
+ \sum_{i=1}^N j(\jx^{i-1}(\mu)) \varepsilon_{m+i}
+ p \sum_{i=1}^N \nu_i \varepsilon_{M+i}.
$$
Note the number $j(\jx^{i-1}(\mu))$ appearing here is simply 
the $i$th part of $\mx(\mu)$ according to Theorem~\ref{xua}.
Finally,  using Theorem~\ref{dt} for the explicit description of
$\Lambda^{++}(M|N,d)$, we deduce that
$r^{-1}(\lambda)$ belongs to $\Lambda^{++}(M|N,d)$
if and only if $j(\mu) \leq \lambda_m$.
Since $j(\mu) = j(\mu + p \nu) = j(t(\lambda))$, this completes the proof.
\end{proof}

\end{document}